\documentclass[10pt, reqno]{amsart}
\usepackage{amssymb}
\usepackage{amsmath,amscd}
\usepackage[initials,nobysame,alphabetic]{amsrefs}
\usepackage{amsmath, amssymb}
\usepackage{amsfonts}
\usepackage{mathrsfs}
\usepackage{mathpazo}
\usepackage[arrow,matrix,curve,cmtip,ps]{xy}
\usepackage[margin=1.3in]{geometry}
\usepackage{amsthm}
\usepackage{float}
\usepackage{hyperref}
\usepackage{graphics}
\usepackage{graphicx}
\usepackage{verbatim}
\usepackage{multirow}
\usepackage{tikz}
\usepackage{enumerate}
\def \qed{\hspace*{\fill} $\Box$\par\medskip}

\def \xtx {X^{t,x}}
\def \xts{X^{t,x}}
\def \vy {\vec y}

\def \usl {\underline}
\def \mi {\mathcal{I}}

\def \txk {(t,x) \in [0,T]\times \R^{k}}

\def \hdd{{\cal H}^{2,d}}
\def \hld{\mathcal{H}^2(\mathcal{L}^2(\lambda))}
\def \cal {\mathcal }
\def \aa {\mathcal{A}^2}
\def \d {\delta}
\def \wt {\mbox{w.r.t}}

\def \ss {\mathcal{S}^2}
\def \ii {i\in \mi}
\def \lb {\label}
\def \tx {(t,x)\in [0,T]\times \R^k}
\def \sol{\underline}
\def \nd {\noindent}
\def \tbf{\textbf}
\def \ed {\end{document}}

\newlength{\oldparindent}
\makeatletter
\newcommand*{\rom}[2]{\expandafter\@slowromancap\romannumeral #1@}
 
\newcommand{\E}{\mathbb{E}}

\newcommand{\eps}{\varepsilon}

\newcommand{\R}{\mathbb{R}}
\newcommand{\N}{\mathbb{N}}

\numberwithin{equation}{section}
\newtheorem{definition}{Definition }[section]
\newtheorem{proposition}[definition]{Proposition }

\newtheorem{theoreme}[definition]
{Theorem }
\newtheorem{corollaire}[definition]
{Corollary }
\newtheorem{remarque}[definition]
{Remark }

\newcommand{\MBFigure}[6]{
$\left. \right.$ \\
\refstepcounter{figure}
\addcontentsline{lof}{figure}{\numberline{\thefigure}{\ignorespaces #5}}
\begin{center}
\begin{minipage}{#1cm}
\centerline{\includegraphics[width=#2cm,angle=#3]{#4}}
\begin{center}
\upshape{F\textsc{ig} \normal
\end{center}
size{\thefigure}. $-$} #5
\end{center}
\label{#6}
\end{minipage}
\end{center}
$\left. \right.$ \\}
\def \g {\gamma}
\def \p {\mathbb{P}}

\def \rw {\rightarrow}

\def \stt {s\in [t,T]}

\def \txk {(t,x) \in [0,T]\times \R^{k}}

\def \fst {\forall s\le T}
\def \TR {[0,T] \times \R}
\def \r{\rho}
\def \xtx {X^{t,x}}
\def \xts{X^{t,x}}
\def \vy {\vec y}
\def \rw {\rightarrow}
\def \fii {\forall \ii}

\def \mi {\mathcal{I}}

\def \txk {(t,x) \in [0,T]\times \R^{k}}

\def \hdd{\mathcal{H}^{2,d}}
\def \g {\gamma}
\def \hld{\mathcal{H}^2(\mathcal{L}^2(\lambda))}
\def \aa {\mathcal{A}^2}
\def \d {\delta}

\def \wt {\mbox{w.r.t}}

\def \ii {i\in \mi}
\def \ki {k\in \mi}
\def \lb {\label}
\def \tx {(t,x)\in [0,T]\times \R^k}
\def \sol{\underline}
\def \nd {\noindent}
\def \tbf{\textbf}
\def \sot {s\in [0,T]}
\def \stt {s\in [t,T]}
\def \tr{\mbox{Tr}}

\def \qed{\hspace*{\fill} $\Box$\par\medskip}

\def \xtx {X^{t,x}}
\def \xts{X^{t,x}}

\def \vy {\vec y}

\def \mi {\mathcal{I}}

\def \txk {(t,x) \in [0,T]\times \R^{k}}

\def \hld{\mathcal{H}^2(\mathcal{L}^2(\lambda))}

\def \t {\tau}

\def \aa {\mathcal{A}^2}
\def \d {\delta}
\def \wt {\mbox{w.r.t}}

\def \ii {i\in \mi}
\def \lb {\label}
\def \tx {(t,x)\in [0,T]\times \R^k}
\def \sol{\underline}
\def \nd {\noindent}
\def \tbf{\textbf}

\def \xtx {X^{t,x}}
\def \xts{X^{t,x}}
\def \vy {\vec y}

\def \mi {\mathcal{I}}

\def \txk {(t,x) \in [0,T]\times \R^{k}}

\def \hld{\mathcal{H}^2(\mathcal{L}^2(\lambda))}
\def \aa {\mathcal{A}^2}
\def \d {\delta}
\def \wt {\mbox{w.r.t}}

\def \ii {i\in \mi}
\def \lb {\label}
\def \tx {(t,x)\in [0,T]\times \R^k}
\def \sol{\underline}
\def \nd {\noindent}
\def \tbf{\textbf}

\title{Viscosity solution of system of integral-partial differential equations with interconnected obstacles of non-local type without Monotonicity Conditions
}
\author{Said Hamad\`ene , Mohamed Mnif, Sarra Neffati 
}
\thanks{\small  LMM, Le Mans Universit\'e, Avenue Olivier Messiaen, 72085 Le Mans, Cedex 9, France, email \href{mailto: said.hamadene@univ-lemans.fr }{ said.hamadene@univ-lemans.fr}  }
\thanks{\small  ENIT, LAMSIN, University of Tunis El Manar, Tunis, Tunisia, email \href{mailto:mohamed.mnif@enit.utm.tn}{mohamed.mnif@enit.utm.tn}   }
\thanks{\small  LMM, Le Mans Universit\'e, Avenue Olivier Messiaen, 72085 Le Mans, Cedex 9, France, email \href{mailto: sarra.neffati@enit.utm.tn  }{ sarra.neffati@enit.utm.tn }  }
\date {September 3, 2024}
\begin{document}
\maketitle
\begin{abstract}
In this paper, we study a system of second order integral-partial differential equations with interconnected obstacles with non-local terms, related to an optimal switching problem in the jump-diffusion model. Getting rid of the monotonicity condition on the generators with respect to the jump component, we construct a  continuous viscosity solution of the system which is unique in the class of functions with polynomial growth. In our study, the main tool is the associated system of reflected backward stochastic differential equations with jumps and interconnected obstacles for which we show the existence of a unique Markovian solution.   
\end{abstract}
	

\section{Introduction}
\makebox[0.75cm][r] Let us consider the following system  of integral-partial differential equations (IPDEs for short) with interconnected obstacles with non-local terms: For any $i \in \mathcal{I}:= \lbrace 1,...,m\rbrace$,
\begin{equation}\label{eqPDE}
\begin{cases}
\min \lbrace u^i(t,x) - \displaystyle \max_{j \in {\mathcal{I}^{-i}}}(u^j(t,x)-g_{ij}(t,x)); -\partial_tu^i(t,x) - \mathcal{L}u^i(t,x) - \mathcal{K}u^i(t,x)\\
 \quad - \bar f_i(t,x,(u^k(t,x))_{k=1,m},(\sigma^\top D_xu^i)(t,x),B_iu^i(t,x))\rbrace = 0, \, \, (t,x) \in [0,T] \times \R^k; \\\\
u^i(T,x) = h_i(x),
\end{cases}
\end{equation}
where $\mathcal{I}^{-i}:= \mathcal{I}-\lbrace i \rbrace$ for any $\ii$ and the operators $\mathcal{L}$, $\mathcal{K}$ and $B_i$ are defined as follows: For any $\ii$,
\begin{equation}\lb{eqintro}
\begin{aligned}
&\mathcal{L}u^i(t,x) :=  b(t,x)^\top D_x u^i(t,x) + \frac{1}{2}\mbox{Tr}[(\sigma\sigma^\top)(t,x)D_{xx}^2 u^i(t,x)],\\[3pt]
&\mathcal{K}u^i(t,x) := \textstyle \int_E (u^i(t,x+\beta(x,e))- u^i(t,x)- \beta(x,e)^\top D_x u^i(t,x))\lambda(de) \, \, \, \mbox{ and }\\[5pt]
& B_iu^i(t,x):= \textstyle \int_E \gamma_i(x,e) (u^i(t,x+\beta(x,e))- u^i(t,x))\lambda(de).
\end{aligned}
\end{equation}

In the above, $D_xu^i$ and $D^2_{xx}u^i$ are the gradient and Hessian matrix of $u^i$ with respect to its second variable $x$, respectively; $(.)^\top$ is the transpose and $\lambda(.)$ is a finite L\'evy measure on $E:= \R^l-\lbrace 0 \rbrace$, i.e., $\lambda(E)<\infty$.

We note that, due to the presence of $B_iu^i$ and $Ku^i$ in equation \eqref{eqPDE}, such an IPDE is called of non-local type. The non-local setting has been studied by several authors (see e.g. \cite{barles1997backward}, \cite{imbert2008}, \cite{hamadene}, \cite{hamadene2016viscosity}, \cite{hamadene2015viscosity}, \cite{niklas2019}). Actually, in \cite{hamadene2015viscosity}, Hamad\`ene-Zhao  have shown that, if for any $\ii$,
 \begin{itemize}
 \item[(i)] $\gamma_i \geqslant 0$,
 \item[(ii)] $ q\in \R \mapsto \bar f_i(t,x,(y_k)_{k=1,m},z,q)$ is non-decreasing, when the other components $(t,x,y,z)$ are fixed,
\end{itemize}
then, there exist functions $(u^i)_{\ii}$ unique continuous viscosity solution of system \eqref{eqPDE} in the class of functions with polynomial growth. Conditions (i)-(ii), which will be referred as the monotonicity conditions, are needed in order to have the comparison result and to treat the  operator $B_iu^i$ which is not well-defined for an arbitrary function $u^i$. The monotonicity conditions (i)-(ii) are usually assumed in the literature on viscosity solutions for equations with a non local term of types \eqref{eqPDE}. Therefore, without assuming the conditions neither on $\gamma_i$ nor on $\bar f_i$, $i=1,...,m$, the problem of existence and uniqueness of the viscosity solution of system \eqref{eqPDE} remains open. To deal with this problem is the main objective of this paper.

A special case of this system of IPDEs with interconnected obstacles occurs in the context of optimal switching control problems when the dynamics of  the state variables are described by a jump diffusion process $(X_s^{t,x})_{s \leq T}$ solution of the following stochastic differential equation:
\begin{equation}\label{SDE}
\begin{cases}
dX_s^{t,x} = b(s, X_s^{t,x}) ds + \sigma(s, X_s^{t,x}) dB_s + \int_E \beta(X_{s{-}}^{t,x}, e) \tilde{\mu}(ds,de),\, \, \, s \in [t,T]; \\
 X_s^{t,x} = x \in \R^k,\, \,    s\leq t,
 \end{cases}
 \end{equation}
where $B := (B_s)_{s \leq T}$ is a $d$-dimensional Brownian motion, $\mu$ an independent Poisson random measure with compensator $ds\lambda(de)$ and $\tilde{\mu}(ds,de):= \mu(ds,de)-ds\lambda(de)$ its compensated random measure.

In this setting, if for any $\ii$, $\bar f_i$ does not depend on $(u^k)_{k=1,m}$, $D_xu^i$ and $B_iu^i$ (see e.g. \cite{hamadene2015systems}), the IPDEs \eqref{eqPDE} reduce to the Hamilton-Jacobi-Bellman system associated with the switching control problem whose value function is defined by: $\forall \ii$ and $\tx$,
\begin{align*}
u^i(t,x)=\sup_{\delta := (\theta_k, \alpha_k)_{k \geq 0}} \E\Big[ \textstyle \int_t^T \bar f^{\d}(s,X_s^{t,x})ds - \sum_{k \geq 1} g_{\alpha_{k-1} \alpha_{k}}(\theta_k, X_{\theta_{k}}^{t,x})\mathbf{1}_{\lbrace \theta_{k} <T \rbrace}
 + h^{\delta}(X_T^{t,x})\Big],
\end{align*}
where :
 
\nd (a) $\delta := (\theta_k, \alpha_k)_{k \geq 0}$ is a strategy of switching in which $(\theta_k)_{k\geq 0}$ is an increasing sequence of stopping times and $(\alpha_k)_{k\geq 0}$ is a sequence of random variables with values in $\mi=\lbrace 1,...,m\rbrace$ $(\theta_0=t$ and $\alpha_0=i$);
 
\nd (b) $ \bar f^{\delta}(s,X_s^{t,x})$ is the instantaneous payoff when the strategy $\d$ is implemented on the system under switching and $h^{\delta}(X_T^{t,x})$ is the terminal payoff ;

\nd (c) $g_{ij}(.)$ is the switching cost function when moving from mode $i$ to mode $j$ $(i, j \in \mi, \, i\neq j)$.

The main tool to tackle system \eqref{eqPDE} is to deal with the following system of reflected backward stochastic differential equations (RBSDEs for short) with jumps and interconnected obstacles: $\forall \ii$ and $s \in [t,T],$
\begin{equation}\label{eqRBSDE}
\begin{split}
\begin{cases}
\vspace{0.3cm}
Y_s^{i,t,x} = h_i(X_T^{t,x})+ \int_s^T \bar f_i(r,X_r^{t,x},(Y_r^{k,t,x})_{k\in \mi},Z_r^{i,t,x}, \int_E V_r^{i,t,x}(e)\gamma_i(X_r^{t,x},e) \lambda(de))dr \\\vspace{0.3cm}
 \qquad \qquad   +K_T^{i,t,x} - K_s^{i,t,x} -\int_s^T Z_r^{i,t,x}dB_r -\int_s^T \int_E V_r^{i,t,x}(e) \tilde{\mu}(dr,de);\\\vspace{0.3cm}
 Y_s^{i,t,x} \geqslant  \displaystyle \max_{j \in \mathcal{I}^{-i}}(Y_s^{j,t,x} -g_{ij}(s,X_s^{t,x}));\\\vspace{0.3cm}
 \textstyle {\int_t^T (Y_s^{i,t,x} -\displaystyle \max_{j\in \mathcal{I}^{-i}}(Y_s^{j,t,x} -g_{ij}(s,X_s^{t,x}))) dK_s^{i,t,x} = 0.}
\end{cases}
\end{split}
\end{equation} 
Note that, without the jump process, the system of RBSDEs with oblique reflection \eqref{eqRBSDE} has been investigated in several papers including
 (\cites{chassagneux2011note,hamadene2010switching,hamadene2013viscosity,hu2010multi}, etc.). With the presence of the jump process,  Hamad\`ene-Zhao in \cite{hamadene2015viscosity}, have proved, under conditions (i)-(ii) on $\gamma_i$ and $\bar f_i$, $i=1,...,m$, the existence and uniqueness of the solution $(Y^{i,t,x}, Z^{i,t,x}, V^{i,t,x}, K^{i,t,x})_{\ii}$ of RBSDEs \eqref{eqRBSDE}. Moreover, they stated the link between this RBSDEs and the IPDEs \eqref{eqPDE} through the Feynman-Kac representation, i.e., for any $\tx$ and $\ii$, 
 \begin{equation}\label{fkc2}\forall s\in [t,T], \,Y^{i,t,x}_s=u^{i}(s,\xtx_s)\mbox{ and }u^{i}(t,x)=Y^{i,t,x}_t.
\end{equation}  
Therefore, in the first part of this paper, the main issue is to deal with RBSDEs \eqref{eqRBSDE} without assuming the two points (i)-(ii) mentioned above. Actually we show that when the measure $\lambda(.)$ is finite, the system of RBSDEs \eqref{eqRBSDE} has a solution which is unique among the Markovian solutions, i.e., which have the representation \eqref{fkc2}. Our method relies mainly on the characterization of the jump part $V^{i,t,x}$ of the RBSDEs \eqref{eqRBSDE} by means of the functions $(u^i)_{i=1,m}$ defined in \eqref{fkc2} and the jump-diffusion process $X^{t,x}$.
 In the second part, we deal with the problem of existence and uniqueness in viscosity sense of the solution of system \eqref{eqPDE}. We show that the functions $(u^i)_{i=1,m}$ defined in \eqref{fkc2}, through the unique solution of 
 \eqref{eqRBSDE}, is the unique viscosity solution of system \eqref{eqPDE}. 
 
The paper is organized as follows. In Section 2, we provide all the necessary notations and assumptions concerning the study of IPDEs \eqref{eqPDE} and related RBSDEs with jumps as well. In Section 3, we study the existence of a solution for system of RBSDEs with jumps \eqref{eqRBSDE} and Feynman-Kac representation \eqref{fkc2}. Actually we introduce an approximating scheme (see system \eqref{recurrence} below) which we show that it converges to the solution of system \eqref{eqRBSDE} when the functions $h_i$ and 
$\bar f_i(t,x,\vec 0,0,0)$, $\ii$, are bounded. On the other hand, the 
Feynman-Kac representation \eqref{fkc2} holds true. Later on, by a weighting technique, we get rid of those latter boundedness conditions on $(h_i)_{\ii}$ and $(\bar f_i(t,x,\vec 0,0,0))_{\ii}$. Finally we show that the Markovian solution of \eqref{eqRBSDE} is unique. At the end of the paper, in Section 4, we prove that the functions $(u^i)_{i=1,m}$, which are connected to $(Y^{i,t,x})_{\ii}$ by \eqref{fkc2}, are the unique viscosity solution of \eqref{eqPDE} in the class of continuous functions with polynomial growth. In the Appendix, we give another definition of the viscosity solution of system \eqref{eqPDE} which is inspired by the work by Hamad\`ene-Morlais in \cite{hamadene2016viscosity}.
 \qed

\section{Preliminaries and notations}
Let $T > 0$ be a given time horizon and $(\Omega, \mathcal{F}, \mathbb{F}:=(\mathcal{F}_t)_{t \leq T},  \mathbb{P})$  be a stochastic basis such that $\mathcal{F}_0$ contains all the $\p$-null sets of $\mathcal{F}$,
 $\mathcal{F}_{t+} = \cap_{\varepsilon > 0 } \mathcal{F}_{t+\varepsilon} = \mathcal{F}_t$, and we suppose that the filtration is generated by the two following mutually independent processes :
 \begin{itemize}
 \item[(i)] a  $d$-dimensional standard Brownian motion $B:= (B_t)_{0 \leq t\leq T}$ ; 
 \item[(ii)] a Poisson random measure $\mu$ on $\R^+ \times E$, where $E:= \R^{l} - \lbrace0\rbrace$ is equipped with its Borel field $\mathcal{B}(E), (l\geqslant 1$ fixed). Let $\nu(dt,de) :=dt \lambda (de)$ be its compensated process. Throughout this paper the measure $\lambda(.)$ is assumed to be finite on $(E, \mathcal{B}(E))$, i.e.,  $\lambda(E)<\infty$.  An example when $l=1$ is $\lambda(de)=(|e|^{-\theta}1_{\{|e|\le 1\}}+|e|^{-\rho}1_{\{|e|\ge 1\}})de$ with $\theta \in (0,1)$ and $\rho>1$. The compensated measure of $\mu$ is denoted by $\tilde \mu$, i.e., $\lbrace \tilde{\mu}([0,t]\times A) := (\mu - \nu)([0,t]\times A)\rbrace_{t\leq T }$ which is a martingale for every $A \in \mathcal{B}(E)$. 
 \end{itemize}
 Let us now introduce the following spaces:
 \begin{itemize}
 \item[a)] $\mathcal{P}$ (resp. $\mathbf{P}$) is the $\sigma$-algebra of $\mathbb{F}$-progressively measurable (resp. $\mathbb{F}$-predictable) sets on $\Omega \times [0,T];$
 \item[b)] $\mathcal{L}^2(\lambda)$ is the space of Borel measurable functions $(\varphi(e))_{e \in E}$ from $E$ into $\R$ such that $\int_E|\varphi(e)|^2 \lambda(de)< \infty$;
 
 \item[c)] $\ss$ is the space of RCLL (right continuous with left limits) $\mathcal{P}$-measurable and $\R$-valued processes $Y:=(Y_s)_{s \leq T}$ such that $\mathbb{E}\big[\displaystyle \sup_{0\leq t\leq T} |Y_s|^2\big] < \infty$;
 \item[d)] $\aa$ is the subspace of $\ss$ of continuous non-decreasing processes $K:= (K_t)_{t\leq T}$ such that $K_0=0$;
 
 \item[e)] $\hdd$ is the space of $\mathcal{P}$-measurable and $\R^{d}$-valued processes $Z:=(Z_s)_{s \leq T}$ such that $\mathbb{E}\left[\int_0^T |Z_s|^2 ds\right] < \infty$; 
 
 \item[f)] $\hld$ is the space of $\mathbf{P}$-measurable and $\mathcal{L}^2(\lambda)$-valued processes $V:=(V_s)_{s \leq T}$ such that\\
  $\mathbb{E}\left[\int_0^T \int_E |V_s(e)|^2 \lambda (de) ds \right]< \infty$. 
 \end{itemize}
 
 For an RCLL process $(\theta_s)_{s\leq T}$, we define for any $s\in (0,T],$
 $\theta_{s-} = \lim_{r\nearrow s}\theta_r$ and $\Delta_s\theta := \theta_s- \theta_{s-}$ is the jump size of $\theta$ at $s$.
 
 Now, for any $(t,x) \in [0,T]\times \R^k$, let  $ (X_s^{t,x})_{s \leq T}$ be the stochastic process solution of the following stochastic differential equation (SDE for short) of diffusion-jump type: 
 \begin{equation}\label{eq1}
 \begin{cases}
 dX_s^{t,x} = b(s, X_s^{t,x}) ds + \sigma(s, X_s^{t,x}) dB_s + \int_E \beta(X_{s{-}}^{t,x}, e) \tilde{\mu}(ds,de), \hspace{0.2cm}  s\in [t, T]; \\
 X_s^{t,x} = x \in \R^k, \hspace{0.3cm} 0\leq s\leq t ,
 \end{cases}
 \end{equation}
 where $b : [0,T]\times \R^k \rightarrow \R^k$ and $\sigma : [0,T]\times \R^k \rightarrow \R^{k\times d}$ are two continuous functions in $(t,x)$ and Lipschitz  $\wt$ $x$, i.e., there exists a positive constant $C$ such that:
 \begin{equation}\label{eq2}
 |b(t,x) - b(t,x^{\prime})| + |\sigma(t,x) - \sigma(t,x^{\prime})|\leq C|x-x^{\prime}|,\hspace{0.2cm} \forall (t,x,x^{\prime}) \in [0,T]\times \R^{k+k}.
 \end{equation}
 The continuity of $b$, $\sigma$ and \eqref{eq2} imply the existence of a constant $C$ such that  \begin{equation}\label{eq3}
 |b(t,x)| + |\sigma(t,x)| \leq C(1+|x|), \hspace{0.2cm} \forall (t,x) \in [0,T]\times \R^{k}.
\end{equation}
The function $\beta : \R^k\times E \rightarrow \R^k$ is measurable and verifies: For some real constant $c$, \begin{align}\label{eq4}
    &|\beta(x,e)|\leq c(1\wedge|e|) \mbox{ and }|\beta(x,e)- \beta(x^{\prime},e)| \leq c|x-x^{\prime}|(1\wedge |e|), \quad  \forall e \in  E \hspace{0.1cm} \mbox{and} \hspace{0.1cm} x, x^{\prime} \in \R^k.
\end{align}
Conditions \eqref{eq2}, \eqref{eq3} and \eqref{eq4} ensure, for any $\txk$, the existence and uniqueness of a  solution of equation \eqref{eq1} (see \cite{fujiwara1985stochastic} for more details). Moreover the following estimate holds true: 
\begin{equation}\label{estimation X} 
\forall p\geqslant 1,\,\,
\mathbb{E}[\sup_{s\leq T} |X_s^{t,x}|^p] \leqslant C(1+|x|^p).
\end{equation}
Next, let us introduce the measurable deterministic functions $(\bar f_i)_{\ii}$, $(h_i)_{\ii}$, $(g_{ij})_{i,j \in \mathcal{I}}$ 
and $(\g_i)_{\ii}$ defined as follows : for any $i,j \in \mathcal{I}$,
\begin{align*}
&a)\,\,\,\bar f_i : (t,x,\vec{y},z,q) \in [0,T] \times \R^{k+m+d+1}  \longmapsto \bar f_i(t,x,\vec y,z,q)\in \R \,\, (\vec y:=(y^1,...,y^m))\,;\\\\
&b)\,\, g_{ij}: \txk \longmapsto g_{ij}(t,x)\in \R^+\,; \\\\&c)\,\,h_i :x\in \R^{k}  \longmapsto h_i(x)\in \R;
\\\\&d)\,\,\g_i :(x,e)\in \R^{k}\times E  \longmapsto \g_i(x,e)\in \R.
\end{align*}
Additionally we assume that they satisfy:
\begin{itemize}
\item[\textbf{(H1)}] For any $i \in \lbrace 1,...,m\rbrace$:
\begin{itemize} 
 \item[(i)]  The function $(t,x) \mapsto \bar f_i(t,x,\vy ,z,q)$  is continuous, uniformly w.r.t. the variables $(\vy, z,q)$.
 \item[(ii)] The function $ \bar f_i$ is Lipschitz continuous w.r.t. the variables $(\vy,z,q)$ uniformly in $(t,x)$, i.e., there exists a positive constant $C_i$ such that for any $(t,x) \in [0,T]\times \R^k,$ 
 $(\vy, z,q)$ and $(\vy_1,z_1,q_1)$ elements of  $\R^{m+d+1}$: 
\begin{equation}
|\bar f_i(t,x,\vy,z,q)- \bar f_i(t,x,\vy_1,z_1,q_1)| \leq C_i( |\vy- \vy_1| + |z - z_1| + |q - q_1|).
\end{equation} 
\item[(iii)] The mapping $(t,x) \mapsto \bar f_i(t,x,0,0,0)$ has polynomial growth in $x$, i.e., there exist two constants $C > 0$ and $p \geqslant 1$ such that for any $(t,x) \in [0,T]\times \R^k$,
\begin{equation}
|\bar f_i(t,x,0,0,0)| \leq C(1+|x|^p).
\end{equation}
\item[(iv)] For any $\ii$ and $ j\in \mathcal{I}^{-i}$,
the mapping $y^j \mapsto \bar f_i(t,x,y^1,...,y^{j-1},y^j,y^{j+1},...,y^m,z,q)$ is non-decreasing whenever the components   $(t,x,y^1,...,y^{j-1},y^{j+1},...,y^m,z,q)$ are fixed.\\
\item[(v)] The functions $ (\gamma_i)_{\ii}$ verify: For any $x,x'$ and $e$,
\begin{equation}\lb{eqgamma}
|\gamma_i(x,e)-\gamma_i(x',e)|\le \bar c^i_\gamma |x-x'|(1\wedge |e|) \text{ and }|\gamma_i(x,e)|\le  c^i_\gamma (1\wedge |e|).
\end{equation}
\end{itemize}
\item[\textbf{(H2)}] $\forall i,j \in \lbrace 1,...,m \rbrace,$ $g_{ii} = 0$ and for $i\neq j,$ $g_{ij}(t,x)$ is non-negative, jointly continuous in $(t,x)$ with polynomial growth and satisfies the following non free loop property :
For any $(t,x) \in [0,T] \times \R^k$, for any sequence of indices $i_1,...,i_k$ such that $i_1 = i_k$ and $card\lbrace i_1,...,i_k \rbrace = k-1$ ($k\ge 3$) we have 
\begin{equation}
g_{i_{1}i_{2}}(t,x) + g_{i_{2}i_{3}}(t,x) + ... + g_{i_{k-1}i_{1}}(t,x)> 0.
\end{equation} 
 \item[\textbf{(H3)}] For $i \in \lbrace 1,...,m\rbrace$, the function $h_i$, which stands for the terminal condition, is continuous with polynomial growth and satisfies the following consistency condition: 
 \begin{equation}\lb{condh}
\forall x \in \R^k,\,\, h_i(x) \geqslant \max_{j\in \mathcal{I}^{-i}}(h_j(x)-g_{ij}(T,x)).
 \end{equation} 
\end{itemize}
\nd We now introduce the following assumptions:
\begin{itemize}
\item[\textbf{(H4)-(i)}] $\forall i \in \mathcal{I}$, $\gamma_i \geqslant 0$;
\item[\textbf{(H4)-(ii)}]  The mapping $
q \in \R \longmapsto \bar f_i(t,x,\vy,z,q)$ is non-decreasing when the other components   $(t,x,\vy,z)$ are fixed.
\qed
\end{itemize}
Next we define the functions $(f_i)_{i=1,...,m}$ on $[0,T] \times \R^{k+m+d} \times \mathcal{L}^2(\lambda),$ as follows: $\fii$, 
 \begin{equation} 
 f_i(t,x,\vy,z,v):= \bar f_i(t,x,\vy,z,\textstyle \int_E v(e)\gamma_i(x,e)\lambda(de)).
\end{equation} 
Note that since $\bar f_i$ is uniformly Lipschitz in $ (\vy, z,q)$ and by (H1)-(v) on $\gamma_i$, the function $f_i$ enjoy the two following properties:
 \begin{itemize}
 \item[(a)] $f_i$ is Lipschitz continuous w.r.t. the variables $(\vy,z,v)$ uniformly in $(t,x)$;
 \item[(b)] The mapping $(t,x) \mapsto f_i(t,x,\vec 0,0,0)=\bar f_i(t,x,\vec{0},0,0)$ is continuous with polynomial growth.
 \end{itemize}

 \begin{remarque}\label{rmqimportante}The condition (H1)-i) is needed, e.g. in \cite{hamadene2016viscosity} or \cite{hamadene2015viscosity} in order to apply Ishii's Lemma to show comparison, of sub- and supersolutions, in the systems considered in those papers and then to deduce uniqueness and continuity of the viscosity solution. However instead of requiring $(H1)-i)$ it is enough to require other sufficient conditions which make comparison of sub. and super-solutions holds
 . So if $\bar f_i$, $\ii$, do not depend on $q$ it is enough to require the following conditions:
\medskip

\noindent (a) For any $\ii$, $\bar f_i$ is jointly continuous in $(t,x,\vy,z)$;

\noindent (b) For any $R>0$, there exists a continuous function $m_R$ from $\R^+$ to $\R^+$ such that $m_R(0)=0$ and for $t$, $|x|\le R$, $|x'|\le R$, $|\vec y|\le R$ and $z$, we have,
\begin{equation}\label{newcondition}
    |\bar f_i(t,x,\vec y,z)-\bar f_i(t,x',\vec y,z)|\le 
    m_R(|x-x'|(1+|z|)).
\end{equation}
One can see e.g. the paper by El-Karoui et al. \cite{kkppq} on this latter condition. 
In the case when $(\bar f_i)_{\ii}$ depend on $q$, similar results exist (one can see e.g. \cite{barles1997backward} for more details). Finally let us notice that the conditions on $\beta$ and $(\gamma_i)_{\ii}$ are not sharp and can be improved since $\lambda(.)$ is finite. However, as the main objective is to get rid of the monotonicity conditions, we prefer 
to make those assumptions which are usually assumed in several papers including \cites{barles1997backward,hamadene2015systems}. 
\end{remarque}
\def \txrt {(t,x)\in [0,T)\times \R^k}
The main objective of this paper is to study the following system of integral-partial differential equations (IPDEs) with interconnected obstacles: for any $i \in \mathcal{I}:=\lbrace 1,...,m\rbrace$, 
\begin{equation}\label{system}
\begin{cases}
\min \lbrace u^i(t,x) - \displaystyle \max_{j \in {\mathcal{I}^{-i}}}(u^j(t,x)-g_{ij}(t,x)); -\partial_tu^i(t,x) - \mathcal{L}u^i(t,x) - \mathcal{K}u^i(t,x)\\
 \qquad - \bar f_i(t,x,(u^k(t,x))_{k=1,m},(\sigma^TD_xu^i)(t,x),B_iu^i(t,x))\rbrace = 0, \txrt; \\
u^i(T,x) = h_i(x),
\end{cases}
\end{equation}
where $\mathcal{L}$ is the second-order local operator
\begin{equation} \label{local-operator}
\mathcal{L}\varphi(t,x) :=  b(t,x)^\top D_x \varphi(t,x) + \frac{1}{2}\mbox{Tr}[(\sigma\sigma^\top)(t,x)D_{xx}^2 \varphi(t,x)].
\end{equation}
The non-local operators $\mathcal{K}$ and 
$B_i$, $\ii$, are defined as follows:
\begin{equation}\label{non-local operators}
\begin{aligned}
&K\varphi(t,x) := \textstyle \int_E (\varphi(t,x+\beta(x,e))- \varphi(t,x)- \beta(x,e)^\top D_x \varphi(t,x))\lambda(de) \, \, \, \mbox{ and }\\[5pt]
& B_i\varphi(t,x):= \textstyle \int_E \gamma_i(x,e) (\varphi(t,x+\beta(x,e))- \varphi(t,x))\lambda(de),
\end{aligned}
\end{equation}
for any $\R$-valued function $\varphi(t,x)$ such that $D_x\varphi(t,x)$ and $D^2_{xx}\varphi(t,x)$ are defined. 
\section{Systems of Reflected BSDEs with Jumps with Oblique Reflection}
The system of IPDEs \eqref{system} is deeply related to the following system of reflected BSDEs with jumps with interconnected obstacles (or oblique reflection) associated with $((\bar{f}_i)_{i\in \mathcal{I}}, (g_{ij})_{i,j \in \mathcal{I}},(h_i)_{i \in \mathcal{I}})$:\\
  $\forall i = 1,...,m$ and $s \in [0,T]$,
 \begin{equation}\label{eqBSDE}
\begin{split}
\begin{cases}
\vspace{0.3cm} Y^{i,t,x}\in \ss,  Z^{i,t,x} \in \hdd, V^{i,t,x} \in \hld, \mbox{ and } K^{i,t,x}\in \aa;\\ \vspace{0.3cm}
Y_s^{i,t,x} = h_i(X_T^{t,x})+ \int_s^T \bar f_i(r,X_r^{t,x},(Y_r^{k,t,x})_{k\in \mi},Z_r^{i,t,x}, \int_E V_r^{i,t,x}(e)\gamma_i(X_r^{t,x},e) \lambda(de))dr \\\vspace{0.3cm}
 \qquad \qquad   +K_T^{i,t,x} - K_s^{i,t,x} -\int_s^T Z_r^{i,t,x}dB_r -\int_s^T \int_E V_r^{i,t,x}(e) \tilde{\mu}(dr,de);\\\vspace{0.3cm}
 Y_s^{i,t,x} \geqslant  \displaystyle \max_{j \in \mathcal{I}^{-i}}(Y_s^{j,t,x} -g_{ij}(s,X_s^{t,x}));\\\vspace{0.3cm}
 \textstyle {\int_0^T (Y_s^{i,t,x} -\displaystyle \max_{j\in \mathcal{I}^{-i}}(Y_s^{j,t,x} -g_{ij}(s,X_s^{t,x}))) dK_s^{i,t,x} = 0.}
\end{cases}
\end{split}
\end{equation}
This system of reflected BSDEs with jumps with interconnected obstacles  \eqref{eqBSDE} has been considered by Hamad\`ene and Zhao in \cite{hamadene2015viscosity} where issues of existence and uniqueness of the solution, and the relationship between the solution of \eqref{eqBSDE} and the one of system \eqref{system}, are considered. Actually in \cite{hamadene2015viscosity}, it is shown:
\begin{theoreme} \label{thmexistence1}(see \cite{hamadene2015viscosity}, pp.1745).\\
 Assume that the deterministic functions $(\bar f_i)_{i\in \mathcal{I}}, (g_{ij})_{i,j \in \mathcal{I}},(h_i)_{i \in \mathcal{I}}$ and $(\gamma_i)_{i \in \mathcal{I}}$  verify Assumptions (H1)-(H3) and (H4). Then, we have:
 \begin{itemize}
 \item[i)] The system \eqref{eqBSDE} has a unique solution $(Y^{i,t,x}, Z^{i,t,x}, V^{i,t,x}, K^{i,t,x})_{i\in \mi}$.
 \item[ii)] There exist deterministic continuous functions  $(u^i)_{\ii}$ of polynomial growth, defined on $[0,T] \times \R^k$, such that:
\begin{equation*}
 \fii, \forall s \in [t,T], \, \, Y_s^{i,t,x} = u^i(s,X_s^{t,x}).
\end{equation*}
\end{itemize} 
\end{theoreme}
In our setting, we also consider the system \eqref{eqBSDE} but without assuming Assumption (H4). 
We then have the following first result which is an intermediary one.
\begin{theoreme} \label{thmexistence2}
Assume that:

\nd (i) The functions $(\bar f_i)_{i\in \mathcal{I}}, (\gamma_i)_{\ii}, (g_{ij})_{i,j \in \mathcal{I}}$ and $(h_i)_{i \in \mathcal{I}}$ verify Assumptions (H1)-(H3).

\nd (ii) There exist a constant $\bar C$ such that, for any $\ii$,
\begin{equation*}
    |h_i(x)|+|\bar f_i(t,x,\vec 0,0,0)| \leq \bar C.
\end{equation*}
Then the system \eqref{eqBSDE} has a solution $(Y^{i,t,x},Z^{i,t,x},V^{i,t,x},K^{i,t,x})_{i\in \mi}$. Moreover there exist  bounded continuous functions $(u^i)_{\ii}$ such that for any $\ii$, $\tx$,
$$Y^{i,t,x}_s=u^i(s,\xtx_s), \,\,\forall s\in [t,T].
$$
\end{theoreme}
\nd {\sol {Proof}}: The proof is divided into four steps.

 \noindent\tbf{Step 1: The iterative construction} 
  
\nd For any $n\geq 0,$ let $(Y^{i,n},Z^{i,n},V^{i,n},K^{i,n})_{i\in \mi}$ be \vspace{0.1cm} the sequence of processes defined recursively as follows:
 \begin{equation*}
   \vspace{0.1cm} (Y^{i,0},Z^{i,0},V^{i,0},K^{i,0}) = (0,0,0,0)\, \, \mbox{ for all } \ii, \, \,  \mbox{ for } n\ge 1 \mbox{ and } s \leq T,
  \end{equation*}
\begin{equation}\label{recurrence}
\begin{split}
\begin{cases}
\vspace{0.2cm} Y^{i,n}\in \ss,  Z^{i,n} \in \hdd, V^{i,n} \in \hld, \mbox{ and } K^{i,n}\in \aa;\\ \vspace{0.2cm}
Y^{i,n}_s = h_i(X_T^{t,x}) +\int_s^T \bar f_i(r,\xtx_r,(Y^{k,n}_r)_{k\in \mi},Z^{i,n}_r, \int_E V^{i,n-1}_r(e)\gamma_i(X_r^{t,x},e)\lambda(de))dr \\ \vspace{0.4cm}
 \qquad \qquad \quad \,  +K_T^{i,n}- K_s^{i,n}-\int_s^T Z_r^{i,n}dB_r -\int_s^T \int_E V_r^{i,n}(e) \tilde{\mu}(dr,de);\\\vspace{0.2cm}
Y_s^{i,n} \geqslant  \displaystyle \max_{j \in \mathcal{I}^{-i}}(Y_s^{j,n} -g_{ij}(s,X_s^{t,x})); \\\vspace{0.2cm}
\textstyle \int_0^T (Y_s^{i,n} -\displaystyle \max_{j\in \mathcal{I}^{-i}}(Y_s^{j,n} -g_{ij}(s,X_s^{t,x}))) dK_s^{i,n} = 0.
\end{cases}
 \end{split}
\end{equation}
First we notice that by Theorem \eqref{thmexistence1}, the solution of this system \eqref{recurrence} exists and is unique. More precisely, for any $\ii,$ the generators $\bar f_i$ do not depend on $ V^{i,n}$, noting that $ V^{i,n-1}$ is already given. The functions $(h_i)_{\ii}$ and $(g_{ij})_{i,j\in \mi}$ satisfy the Assumptions (H2)-(H3) as well. Next, since the setting is  Markovian and using an induction argument on $n$, there exist deterministic continuous functions of polynomial growth  $u^{i,n} : [0,T] \times \R^k \rightarrow \R$, such that for any $s \in [t,T]$:
 \begin{align}\label{rep1}
     &\mbox{(a) }\,Y_s^{i,n} := u^{i,n}(s,\xtx_s)\, \mbox{and }\nonumber\\[4pt]
     & \mbox{(b) }\, V_s^{i,n}(e) := u^{i, n}(s,\xtx_{s{-}}+ \beta(\xtx_{s{-}},e))- u^{i,n}(s,\xtx_{s{-}}),\,\,ds\otimes d\mathbb{P} \otimes d\lambda \mbox{ on } [t,T] \times \Omega \times E.
 \end{align}
 Indeed, for $n=0$, the representations (a), (b) are valid with $u^{i,0}=0$, $\ii$. Assume now that they are satisfied for some $n-1,$ with $n \geq 1$. Then $(Y^{i,n},Z^{i,n},V^{i,n},K^{i,n})$ verifies: for any $s \in [t,T]$ and $\ii$,
\begin{align*}
\begin{cases}
 &Y_s^{i,n} = h_i(X_T^{t,x}) +\int_s^T \bar f_i(r,X_r^{t,x},(Y_r^{k,n})_{k \in \mi},Z_r^{i,n},\int_E \{u^{i, n-1}(r,\xtx_{r-}+ \beta(\xtx_{r-},e))\\[7pt]
 &\qquad \qquad - u^{i,n-1}(r,\xtx_{r-})\}\gamma_i(\xtx_r,e)\lambda(de))dr+K_T^{i,n}- K_s^{i,n}
 - \int_s^T Z_r^{i,n}dB_r -\int_s^T \int_E V_r^{i,n}(e) \tilde{\mu}(dr,de);\\[7pt]
&Y_s^{i,n} \geqslant  \displaystyle \max_{j \in \mathcal{I}^{-i}}(Y_s^{j,n} -g_{ij}(s,X_s^{t,x})); \\[7pt]
 &\int_t^T (Y_s^{i,n} -\displaystyle \max_{j\in \mathcal{I}^{-i}}(Y_s^{j,n} -g_{ij}(s,X_s^{t,x}))) dK_s^{i,n} = 0.
\end{cases}
\end{align*} 
Hence, by Theorem \eqref{thmexistence1}, we deduce the existence of $u^{i,n}$, continuous and of polynomial growth, such that $Y_s^{i,n} = u^{i,n}(s,\xtx_s)$, $s\in [t,T]$ .
Finally as the measure $\lambda$ is finite, i.e., $\lambda(E) < \infty$, then we have the following relationship between the process $(V^{i,n})_{\ii}$ and the deterministic functions $(u^{i,n})_{\ii}$ (see \cite{hamadene2016viscosity}, Proposition 3.3): 
 \begin{equation} \label{repV}
\vspace{0.1cm} V_s^{i,n}(e) = u^{i, n}(s,\xtx_{s{-}}+ \beta(\xtx_{s{-}},e))- u^{i,n}(s,\xtx_{s{-}}),\,\, ds\otimes d\mathbb{P} \otimes d\lambda \mbox{ on } [t,T] \times \Omega \times E.
\end{equation} 
Thus, the two representations (a) and (b) hold true for any $n \geq 0$. Here let us point out that the functions $
(\bar F_i(s,x,\vec y,z):=\bar f_i(s,x,\vec y,z,\int_E \{u^{i, n-1}(s,x+ \beta(x,e))- u^{i,n-1}(s,x)\}\gamma_i(x,e)\lambda(de)))_{\ii}$ verify the assumption (H1).
 \begin{remarque}\label{surot}
 For $s\in [0,t]$, $\xtx_s=x$ and 
 $Y_t^{i,n}=u^{i,n}(t,x)$, therefore in considering the declination of system \eqref{recurrence} on the time interval $[0,t]$, we can easily show by induction that 
 $Z^{i,n}_s1_{[s\le t]}=0, \,\,ds\otimes d\mathbb{P}-a.e$ and 
 $V^{i,n}_s(e)1_{[s\le t]}=0,\,\,  ds\otimes d\mathbb{P}\otimes d\lambda-a.e.$ since the data are continuous and deterministic on $[0,t]$.
 \end{remarque}
\noindent \tbf{Step 2: Switching representation}

 In this step, we represent $Y^{i,n}$ as the value of  an optimal switching problem. Indeed, let $\delta:= (\theta_k, \alpha_k)_{k \geq 0}$ be an admissible strategy of switching, i.e., $(\theta_k)_{k \geq 0}$ is an increasing sequence of stopping times with values in $[0,T]$ such that $\mathbb{P}[\theta_k < T, \forall k \geq 0] = 0$ and $ \forall k \geq 0$, $\alpha_k$ is a random variable $\mathcal{F}_{\theta_{k}}$-measurable with values in $\mathcal{I}.$

 Next, with the admissible strategy $\delta:= (\theta_k, \alpha_k)_{k \geq 0}$ is associated a switching cost process $(A_s^{\delta})_{s \leq T}$ defined by:
 \begin{equation}\label{defcoutswitch}
 A^{\delta}_s := \displaystyle \sum_{k \geq 1} g_{\alpha_{k-1} \alpha_{k}}(\theta_k, X_{\theta_{k}}^{t,x})\mathbf{1}_{\lbrace \theta_{k} \leq s \rbrace} \mbox{ for }s < T,\mbox{ and } A^{\delta}_T= \lim_{s\rightarrow T}A^{\delta}_s.
 \end{equation}
Note that the process $(A^{\delta}_s)_{s \leq T}$ is non-decreasing and RCLL. Finally, for any fixed $s\leq T$ and $i \in \mathcal{I}$, let us denote by $\mathcal{A}_s^i$ the following set of admissible strategies:
\begin{align*}
\mathcal{A}_s^i:= \lbrace\delta:= (&\theta_k, \alpha_k)_{k \geq 0} \ \mbox{admissible strategy such that}\ \theta_0=s, \alpha_0=i \\
&\mbox{and}\ \mathbb{E}[(A^{\delta}_T)^2] < \infty \rbrace.
\end{align*}
 Now, let  $\delta:= (\theta_k, \alpha_k)_{k \geq 0} \in \mathcal{A}_s^i$ and let us define the triplet of adapted processes $(P_\t^{n,\d},N_\t^{n,\d},Q_\t^{n,\d})_{\t \leq T}$ as follows: $\forall \t \leq T,$
 \begin{equation}\label{eq8}\left\{\begin{array}{l}
 P^{n,\d} \mbox{ is RCLL and } \E[\sup_{\t\le T}|P^{n,\d}_\t|^2]<\infty \,; 
 N^{n,\delta}\in \hdd \mbox{ and }  Q^{n,\delta}\in \hld ;\\\\ 
P_\t^{n,\delta} = h^{\delta}(X_T^{t,x})  - A_T^{\delta} + A_\t^{\delta} -\textstyle \int_\t^T N_r^{n,\d}dB_r - \textstyle\int_\t^T \int_E Q_r^{n,\d}(e)\tilde{\mu}(dr,de) \\\\\ \textstyle +\int_\t^T \bar f^{\delta}(r,X_r^{t,x},(Y_r^{k,n})_{k \in \mi},N_r^{n,\delta},
\underbrace{\textstyle\int_E \lbrace u^{\d, n-1}(r,\xtx_{r-}+ \beta(\xtx_{r-},e)) - u^{\d,n-1}(r,\xtx_{r-})\rbrace 1_{\{r\ge s\}} \gamma^{\d}(r,X_r^{t,x},e)  \lambda(de)}_{\Sigma_r^{\d,n-1}})dr ; \end{array}\right.
 \end{equation}
where $h^\d(x)=\sum_{k\ge 0}h^{\alpha_k}(x)1_{\{\theta_k< T\le \theta_{k+1}\}}$, $\gamma^\d(r,\xtx_r,e)=\sum_{k\ge 0}\gamma_{\alpha_k}(\xtx_r,e)1_{\{\theta_k<r\le \theta_{k+1}\}}$,\\
$u^{\d,n-1}(.)=\sum_{k\ge 0}u^{\alpha_k,n-1}(.)1_{\{\theta_k<r\le \theta_{k+1}\}}$
and finally 
$\bar f^{\d}(.)=\sum_{k\ge 0}\bar f_{\alpha_k}(.)1_{\{\theta_k<r\le \theta_{k+1}\}}$. Those series contain only a finite many terms as $\d$ is admissible and then $\mathbb{P}[\theta_n<T,\forall n\geq 0]=0$.
Note that, in \eqref{eq8}, the generators $\bar f^{\d}$ do not depend neither on $P^{n,\d}$ nor on $Q^{n,\d} \in \hld.$ Now, by a change of variables, the existence of $(P^{n,\d} - A^{\d},N^{n,\d}, Q^{n,\d})$ stems from the standard existence result of solutions of BSDEs with jumps by Tang-Li (\cite{tang1994necessary}, pp. 1455) since its generator $ z\mapsto \bar f^{\delta}(r,X_r^{t,x},(Y_r^{k,n})_{k \in \mi},z,\Sigma_r^{\d,n-1})$ is Lipschitz $\wt$ $z$ and $A_T^\d$ is square integrable.

Next, let us consider the following system of RBSDEs: $ \forall \ii$ and $s\leq T$
\begin{equation}\label{BSDE}
\begin{cases}
 \underbar Y^{i,n}\in \ss,  \underbar Z^{i,n} \in \hdd, \underbar V^{i,n} \in \hld, \mbox{ and } \underbar K^{i,n}\in \aa;\\[7pt]
\underbar Y_s^{i,n} = h_i(X_T^{t,x})+ \int_s^T \bar f_i(r,X_r^{t,x},(Y_r^{k,n})_{k\in \mi},\underbar Z_r^{i,n},\int_E V^{i,n-1}_r(e)\gamma_{i}(X_r^{t,x},e)  \lambda(de))dr \\[7pt]\qquad \qquad \qquad + \underbar K_T^{i,n} - \underbar K_s^{i,n} -\int_s^T \underbar Z_r^{i,n}dB_r -\int_s^T \int_E \underbar V_r^{i,n}(e) \tilde{\mu}(dr,de);\\[7pt]
 \underbar Y_s^{i,n} \geqslant  \displaystyle \max_{j \in \mathcal{I}^{-i}}(\underbar Y_s^{j,n} -g_{ij}(s,\xts_s));\\[7pt]
 \textstyle {\int_0^T (\underbar Y_s^{i,n} -\displaystyle \max_{j\in \mathcal{I}^{-i}}(\underbar Y_s^{j,n} -g_{ij}(s,\xtx_s))) d\underbar K_s^{i,n} = 0}
\end{cases}
\end{equation}
  whose solution exists and is unique by Theorem \ref{thmexistence1}. Therefore,  we have the following representation of $\underbar Y^{i,n}$  (see e.g. \cite{hamadene2015systems} for more details on this representation):
 $$\underbar Y_s^{i,n}=\mbox{esssup}_{\d \in \mathcal{A}_s^{i}}(P_s^{n,\d} - A_s^{\d}). $$
But $(Y^{i,n}, Z^{i,n}, V^{i,n},K^{i,n})_{\ii}$ is also solution of \eqref{BSDE}, then by uniqueness one deduces that                     
\begin{equation}\label{eq9}
\forall s\le T,\,\, \underbar Y_s^{i,n} = Y_s^{i,n} = \displaystyle \mbox{esssup}_{\d \in \mathcal{A}_s^i}(P_s^{n,\d} - A_s^{\d}) = (P_s^{n,\delta^{*}} - A_s^{\delta^{*}}),
 \end{equation}
for some $\delta^{*} \in \mathcal{A}_s^i$, which means that $\d^*$ is an optimal strategy of the switching control problem.

\noindent \tbf{Step 3: Convergence of} $\mathbf{(u^{i,n})_{n \geq 0}}$

We now adapt the argument already used in \cites{chassagneux2011note,hamadene2016viscosity,hamadene2015systems} to justify the convergence of the sequence $((u^{i,n})_{\ii})_{n\geq 0}$. For this, let us set for $\ii$ and $n,\, p \geq 1$,
\begin{align*}
F_i^{n,p}( s,\omega,\xtx_s,  z) &:=\bar f_i(s,\xtx_s, (Y_s^{k,n})_{k\in \mi}, z, \textstyle\int_E V_s^{i,n-1}(e) \gamma_{i}(X_s^{t,x},e)  \lambda(de))\\[3pt] & \quad \vee\bar f_i(s,\xtx_s, (Y_s^{k,p})_{k\in \mi}, z, \textstyle\int_E  V_s^{i,p-1}(e) \gamma_{i}(X_s^{t,x},e)  \lambda(de)).
\end{align*}
Next,  let us consider the solution, denoted by  $(\hat Y^{i,n,p}, \hat Z^{i,n,p}, \hat V^{i,n,p},\hat K^{i,n,p})_{\ii}$, of the obliquely reflected BSDEs with jumps associated with $((F_i^{n,p})_{i\in \mathcal{I}}, (g_{ij})_{i,j \in \mathcal{I}}$,  
$(h_i)_{i \in \mathcal{I}})$, which exists and is unique according to Theorem \ref{thmexistence1}. Moreover, as in \eqref{eq9}, we have: $\forall s \leq T$,
\begin{equation}\label{eq10}
 \hat {Y}_s^{i,n,p} = \displaystyle \mbox{esssup}_{\delta \in \mathcal{A}_s^i}(\hat {P}_s^{n,p,\d} - A_s^{\d})=(\hat {P}_s^{n,p,{\tilde \d}^*} - A_s^{{\tilde \delta}^*}),
\end{equation}
where $(\hat{P}^{n,p,\d}, \hat N^{n,p,\d}, \hat Q^{n,p,\d})$ is the solution of the BSDE \eqref{eq8} with generator $F^{\d,n,p}(...)$ which is defined through $F^{n,p}_i(.)$ as the definition of $\bar f^{\d}(.)$ of \eqref{eq8}. Then by the comparison result (see Proposition 4.2 in \cite{hamadene2015viscosity}), between the solutions $Y^{i,n}$ and $\hat{Y}^{i,n,p}$, and $Y^{i,p}$ and $\hat{Y}^{i,n,p}$ (this is possible since the generators of the systems do not depend on the jump parts),  one deduces that: $\forall \ii$ and $s\leq T$,
  \begin{equation*}
  Y_s^{i,n} \leq \hat {Y}_s^{i,n,p}\, \, \,   \mbox{ and } \, \, \,   Y_s^{i,p} \leq \hat {Y}_s^{i,n,p}. 
  \end{equation*}
   This combined with \eqref{eq9} and \eqref{eq10}, lead to:
   \begin{equation*}
       P_s^{n,{\tilde \delta}^*} - A_s^{{\tilde \delta}^*}\leq  Y_s^{i,n} \leq \underbrace{\hat {P}_s^{n,p,{\tilde \d}^*} - A_s^{{\tilde \delta}^*} }_{\hat Y^{i,n,p}_s}\mbox{ and } \,  P_s^{p,{\tilde \delta}^*} - A_s^{{\tilde \delta}^*}\leq  Y_s^{i,p} \leq \underbrace{\hat {P}_s^{n,p,{\tilde \d}^*} - A_s^{{\tilde \delta}^*} }_{\hat Y^{i,n,p}_s}.
   \end{equation*} It implies that 
\begin{equation}\label{eq12} \,\,
\vspace{0.2cm}|Y_s^{i,n} - Y_s^{i,p}|^2 \leq 2 \{|\hat {P}_s^{n,p,{\tilde \delta}^*} - P_s^{n,{\tilde \delta}^*}|^2 + |\hat {P}_s^{n,p,{\tilde \delta}^*} - P_s^{p,{\tilde \delta}^*}|^2\}.
\end{equation}
Since both terms on the right-hand side of \eqref{eq12} are treated similarly, we focus only on the first one. Applying It\^o's formula with $e^{\alpha s}|\hat {P}_s^{n,p,{\tilde \delta}^*} - P_s^{n,{\tilde \delta}^*}|^2$  $(\alpha>0)$, yields: $\forall s \leq T$,
\begin{align}\label{3.10}
  & e^{\alpha s}|\hat {P}_s^{n,p,{\tilde \delta}^*} - P_s^{n,{\tilde \delta}^*}|^2  + \textstyle\int_s^T  e^{\alpha r}|\hat {N}_r^{n,p,{\tilde \delta}^*} - N_r^{n,{\tilde \delta}^*}|^2dr
   +\sum_{s<r\leq T} e^{\alpha r}  \Delta_r(\hat {P}^{{n,p,{\tilde \delta}^*}} - P^{n,{\tilde \delta}^*})^2\nonumber \\[6pt]
  & = -\alpha\textstyle \int_s^T e^{\alpha r}|\hat {P}_r^{n,p,{\tilde \delta}^*} - P_r^{n,{\tilde \delta}^*}|^2dr + 2 \int_s^T  e^{\alpha r}(\hat {P}_r^{n,p,{\tilde \delta}^*} - P_r^{n,{\tilde \delta}^*}) \big\lbrace F^{{\tilde \delta^*},n,p}(r,X_r^{t,x},\hat {N}_r^{n,p,{\tilde \delta}^*})\nonumber \\[6pt]
  & \quad - \bar f^{{\tilde \delta}^*}(r,\xtx_r, (Y_r^{k,n})_{k\in \mi},  {N}_r^{n,{\tilde \delta}^*}, \textstyle \int_E V_r^{{\tilde \delta^*},n-1}(e) \gamma^{{\tilde \delta}^*}(r,X_r^{t,x},e)  \lambda(de))\big\rbrace dr\nonumber \\[6pt]
  & \quad -2 \textstyle\int_s^T e^{\alpha r} (\hat {P}_r^{n,p,{\tilde \delta}^*} - P_r^{n,{\tilde \delta}^*}) (\hat {N}_r^{n,p,{\tilde \delta}^*}-  N_r^{n,{\tilde \delta}^*})dB_r-2 \textstyle\int_s^T\int_E e^{\alpha r} (\hat {P}_{r-}^{n,p,{\tilde \delta}^*} - P_{r-}^{n,{\tilde \delta}^*})(\hat {Q}_r^{n,p,{\tilde \delta}^*}(e) -  Q_r^{n,{\tilde \delta}^*}(e))\tilde{\mu}(de,dr).
\end{align}
Observe that  the inequality  $|x\vee y-x|\le |x-y|$, $\forall x,y\in \R$, combined with the Lipschitz property of $\bar f^{{\tilde \delta}^*}(.)$  lead to:
\begin{align*}
     &|F^{{\tilde \delta^*},n,p}(r,X_r^{t,x},\hat {N}_r^{n,p,{\tilde \delta}^*})- \bar f^{{\tilde \delta}^*}(r,\xtx_r, (Y_r^{k,n})_{k\in \mi},  {N}_r^{n,{\tilde \delta}^*}, \textstyle \int_E V_r^{{\tilde \delta^*},n-1}(e) \gamma^{{\tilde \delta^*}}(r,X_r^{t,x},e)  \lambda(de))|\\[6pt]
     & \qquad\qquad \leq C \big\lbrace |(Y_r^{k,n})_{k\in \mi}-(Y_r^{k,p})_{k\in \mi}| +|\hat {N}_r^{n,p,{\tilde \delta}^*}- {N}_r^{n,{\tilde \delta}^*}|+|\textstyle\int_E\{ V_r^{\tilde \delta^*,n-1}(e)- V_r^{{\tilde \delta}^*,p-1}(e)\}\gamma^{{\tilde \delta}^*}(r,X_r^{t,x},e)  \lambda(de)|\big\rbrace.
\end{align*}\\
Going back to \eqref{3.10}, taking expectation and using the inequality $2|ab|\leq \eps|a|^2 + \frac{1}{\eps}|b|^2 (\eps>0)$, we obtain: 
\begin{align*}
     &\E \Big[ \textstyle e^{\alpha s}|\hat {P}_s^{n,p,{\tilde \delta}^*} - P_s^{n,{\tilde \delta}^*}|^2 + \int_s^T  e^{\alpha r}|\hat {N}_r^{n,p,{\tilde \delta}^*} - N_r^{n,{\tilde \delta}^*}|^2dr+\textstyle \sum_{s<r\leq T} e^{\alpha r}   \Delta_r(\hat {P}^{{n,p,{\tilde \delta}^*}} - P^{n,{\tilde \delta}^*})^2\Big]\\[6pt]
     &\qquad \leq (-\alpha + 3 \eps) \E\Big[\textstyle \int_s^T e^{\alpha r}|\hat {P}_r^{n,p,{\tilde \delta}^*} - P_r^{n,{\tilde \delta}^*}|^2dr\Big]\\[6pt]
     &\qquad + C^2\eps^{-1}\Big\lbrace  \E\Big[ \textstyle \int_s^T  e^{\alpha r}|(Y_r^{k,n})_{k\in \mi}-(Y_r^{k,p})_{k\in \mi}|^2dr \Big]  +\E \Big[ \textstyle\int_s^T e^{\alpha r}|\hat {N}_r^{n,p,{\tilde \delta}^*}- {N}_r^{n,{\tilde \delta}^*}|^2 dr\Big]\\[6pt]
     & \qquad + \E\Big[\int_s^T e^{\alpha r} \big(\int_E |\{V_r^{{\tilde \delta}^*,n-1}(e)-V_r^{{\tilde \delta}^*,p-1}(e)\} \gamma^{{\tilde \delta}^*}(r,\xtx_r,e)| \lambda(de)\big)^2dr\Big]\Big \}.
\end{align*}
If we choose $\alpha = \alpha_0 = 3 \eps$ and $\eps> C^2$, we get : $\forall s \leq T$,
\begin{align*}
    &\E \Big[ e^{\alpha_0 s}|\hat {P}_s^{n,p,{\tilde \delta}^*} - P_s^{n,{\tilde \delta}^*}|^2 \Big] \leq 3C^2\alpha_0^{-1}\Big\lbrace  \E\Big[ \textstyle \int_s^T e^{\alpha_0 r} |(Y_r^{k,n})_{k\in \mi}-(Y_r^{k,p})_{k\in \mi}|^2dr \Big]\\[4pt]
    & \qquad + \E\Big[\textstyle \int_s^T e^{\alpha_0 r} \big(\int_E \sum_{k=1,m}|\{V_r^{k,n-1}(e)-V_r^{k,p-1}(e)\} \gamma_k(\xtx_r,e)| \lambda(de)\big)^2dr\Big]\Big \}.
\end{align*}
The same reasoning leads to the same estimate for  $e^{\alpha s}|\hat {P}_s^{n,p,{\tilde \delta}^*} - P_s^{p,{\tilde \delta}^*}|^2$. Therefore, we deduce from \eqref{eq12} that:
\begin{align}\lb{eq9x}
  \nonumber &\E \Big[  e^{\alpha_0 s}|Y_s^{i,n}- Y_s^{i,p}|^2\Big]  \leq 6C^2\alpha_0^{-1}\Big\lbrace  \E\Big[ \textstyle \int_s^T  e^{\alpha_0 r}|(Y_r^{k,n})_{k\in \mi}-(Y_r^{k,p})_{k\in \mi}|^2dr \Big]\\[6pt]
    & \qquad+ \E\Big[\textstyle \int_s^T e^{\alpha_0 r} \big(\int_E \sum_{k=1,m}|\{V_r^{k,n-1}(e)-V_r^{k,p-1}(e)\} \gamma_k(\xtx_r,e)| \lambda(de)\big)^2dr\Big]\Big \}.
\end{align}
Then, by summing over $i \in \mathcal{I}$, there exists a constant $\kappa$ such that: $\forall s \leq T$,
\begin{align*}
     &\E \Big[  e^{\alpha_0 s}|(Y_s^{k,n})_{k\in \mathcal{I}}- (Y_s^{k,p})_{k\in \mathcal{I}}|^2\Big]  \leq \kappa\Big\lbrace  \E\Big[ \textstyle \int_s^T e^{\alpha_0 r} |(Y_r^{k,n})_{k\in \mi}-(Y_r^{k,p})_{k\in \mi}|^2dr \Big]\\[6pt]
    & \qquad +\E\Big[\textstyle \int_s^T e^{\alpha_0 r} \big(\int_E \sum_{k=1,m}|\{V_r^{k,n-1}(e)-V_r^{k,p-1}(e)\} \gamma_k(\xtx_r,e)| \lambda(de)\big)^2dr\Big]\Big \}.
\end{align*}
Finally by using Gronwall's inequality one can find a constant $\kappa_1$ such that: $\forall s \leq T$,
\begin{align*}
     &\E \Big[  e^{\alpha_0 s}|(Y_s^{k,n})_{k\in \mathcal{I}}- (Y_s^{k,p})_{k\in \mathcal{I}}|^2\Big] \\[6pt]
     & \qquad \leq \kappa_1 \E\Big[\textstyle \int_s^T e^{\alpha_0 r} \big(\int_E \sum_{k=1,m}|\{V_r^{k,n-1}(e)-V_r^{k,p-1}(e)\} \gamma_k(\xtx_r,e)| \lambda(de)\big)^2dr\Big].
\end{align*}
Taking $s=t$ and considering \eqref{rep1}-((a),(b)), we obtain : for any $\ii$,
\begin{align*}
 |u^{i,n}(t,x)- u^{i,p}(t,x)|^2 &\le \sum_{k=1,m}|u^{k,n}(t,x)- u^{k,p}(t,x)|^2   \nonumber \\[6pt]
& \leq \kappa_1  \E\Big[\textstyle \int_t^T e^{\alpha_0 (r-t)} \big(\int_E \sum_{k=1,m} |\{ u^{k, n-1}(r,\xtx_{r-}+ \beta(\xtx_{r-},e)) - u^{k,n-1}(r,\xtx_{r-})\nonumber\\[6pt]
& \quad - ( u^{k, p-1}(r,\xtx_{r-}+ \beta(\xtx_{r-},e))-u^{k,p-1}(r,\xtx_{r-}))\}\gamma_k(\xtx_r,e) |\lambda(de)\big)^2dr\Big].
\end{align*}
 Next, using Cauchy-Schwarz inequality, \eqref{eqgamma} and the inequality $|a+b|^2\leq 2(|a|^2+|b|^2)$, we get:
\begin{align}\label{estimationu}
 &|u^{i,n}(t,x)- u^{i,p}(t,x)|^2 \nonumber\\[4pt]
 &\quad \leq \kappa_1
\E \big[  \textstyle \int_t^T e^{\alpha_0 (r-t)} \big(\int_E \{\sum_{k=1,m}\big|\gamma_k(\xtx_r,e)\big|^2 \}\lambda(de) \big)\nonumber\\[6pt]
& \quad\qquad\qquad\ \times \big(\textstyle\int_E \{\textstyle \sum_{k=1,m}\big| u^{k, n-1}(r,\xtx_{r{-}}+ \beta(\xtx_{r{-}},e))- u^{k,n-1}(r,\xtx_{r{-}})\nonumber\\[6pt]
& \quad\qquad\qquad-  u^{k, p-1}(r,\xtx_{r{-}} + \beta(\xtx_{r{-}},e))+u^{k,p-1}(r,\xtx_{r{-}}) \big|^2\} \lambda(de) \big)dr\big]\nonumber \\[6pt]
&\quad\leq C  \E \big[ \textstyle \int_t^T e^{\alpha_0 (r-t)}\big(\int_E \{\sum_{k=1,m}\big| u^{k, n-1}(r,\xtx_{r{-}}+ \beta(\xtx_{r{-}},e)) - u^{k,n-1}(r,\xtx_{r{-}})\nonumber \\[6pt]
& \quad\qquad\qquad-  u^{k, p-1}(r,\xtx_{r{-}}+ \beta(\xtx_{r{-}},e))+u^{k,p-1}(r,\xtx_{r{-}}) \big|^2\} \lambda(de) \big)dr\big]\nonumber\\[6pt]
& \quad\leq  2C  \E \big[ \textstyle \int_t^T e^{\alpha_0 (r-t)}\int_E \sum_{k=1,m}\{\big| (u^{k, n-1}-u^{k, p-1}) (r,\xtx_{r{-}}+ \beta(\xtx_{r{-}},e))\big|^2 \nonumber\\[7pt]
&\qquad \qquad\quad+ \big|(u^{k,n-1}- u^{k,p-1})(r,\xtx_{r{-}})\big|^2\} \lambda(de)dr\big],
\end{align}
for some constant $C$ (which may change from line to line).\\ 
Now, in order to take the supremum on the inequality \eqref{estimationu}, we need to show the boundedness of $(u^{i,n})_{\ii}$.  For this, 
let $(\bar Y, \bar Z)$ be the solution of the following standard BSDE: for any $s\leq T$, 
\begin{equation*}
    \begin{cases}
    \bar Y \in \ss,\, \bar Z \in \hdd;\\
    \bar Y_s = \bar C + \int_s^T \big\lbrace \bar C+ m\,C^y_f\,\bar Y_r + C^z_f\,|\bar Z_r| +2\theta \bar Y_r\big\rbrace dr - \int_s^T \bar Z_rdB_r;
    \end{cases}
\end{equation*}
where $C_f^y$ , $C_f^z$ and $C_f^v$ are the maximum of the Lipschitz constants of the $\bar f_i(.)'s$ w.r.t. $\vy$, $z$ and $v$ respectively, and  
\begin{equation*}
\begin{array}{c}
    \theta = C^v_f\,C_\gamma\int_E(1\wedge|e|)\lambda(de).
    \end{array}
\end{equation*}
where $C_\g=\max_{\ii}c_\g^i$ ($c_\g^i$ is defined in \eqref{eqgamma}). The solution of this BSDE  exists and is unique by Pardoux and Peng's result \cite{Pardoux-Peng1}. Then, there exists a constant $\underbar C$ such that $|\bar Y|\leq \underbar C$. Finally, noting that $\bar Y$ is deterministic and $\bar Z =0$.\\
 \textbf{}\\
 Now, recall that $((Y^{i,n},Z^{i,n},V^{i,n},K^{i,n})_{i\in \mi})_{n \geq 0}$ verify: 
 \begin{equation*}
    (Y^{i,0},Z^{i,0},V^{i,0},K^{i,0}) = (0,0,0,0)  \mbox{ and for } n\geq 1,
  \end{equation*}
\begin{equation}\label{recc2}
\begin{split}
\begin{cases}
\vspace{0.2cm}
Y^{i,n}_s = h_i(X_T^{t,x}) +\int_s^T \bar f_i(r,\xtx_r,(Y^{k,n}_r)_{k\in \mi},Z^{i,n}_r, \int_E V^{i,n-1}_r(e)\gamma_i(X_r^{t,x},e)\lambda(de))dr\, \\ \vspace{0.2cm}
 \qquad \qquad   +K_T^{i,n}- K_s^{i,n}-\int_s^T Z_r^{i,n}dB_r -\int_s^T \int_E V_r^{i,n}(e) \tilde{\mu}(dr,de), \quad s\leq T;\\\vspace{0.2cm}
Y_s^{i,n} \geqslant  \displaystyle \max_{j \in \mathcal{I}^{-i}}(Y_s^{j,n} -g_{ij}(s,X_s^{t,x})),\quad s\leq T; \\\vspace{0.2cm}
\textstyle \int_0^T (Y_s^{i,n} -\displaystyle \max_{j\in \mathcal{I}^{-i}}(Y_s^{j,n} -g_{ij}(s,X_s^{t,x}))) dK_s^{i,n} = 0,
\end{cases}
 \end{split}
\end{equation}
with, $\forall s\in [t,T]\,$, $Y_s^{i,n}=u^{i,n}(s,X_s^{t,x})$. Then, by an induction argument on $n$, 
we have that: $\forall n\geq 1$ and $\ii$, 
\begin{equation}\label{induction}
 \forall (t,x)\in [0,T]\times \R^k,\quad  |u^{i,n}(t,x)|\leq \bar Y_t.
\end{equation}
Indeed, for $n=1$, we have: 
\begin{equation}\label{3.14}
\begin{split}
\begin{cases}
\vspace{0.2cm}
Y^{i,1}_s = h_i(X_T^{t,x}) +\int_s^T \bar f_i(r,\xtx_r,(Y^{k,1}_r)_{k\in \mi},Z^{i,1}_r,0)dr +K_T^{i,1}- K_s^{i,1}\\ \vspace{0.2cm}
 \qquad \qquad  \, -\int_s^T Z_r^{i,1}dB_r -\int_s^T \int_E V_r^{i,1}(e) \tilde{\mu}(dr,de), \quad s\leq T;\\\vspace{0.2cm}
Y_s^{i,1} \geqslant  \displaystyle \max_{j \in \mathcal{I}^{-i}}(Y_s^{j,1} -g_{ij}(s,X_s^{t,x})),\quad s\leq T; \\\vspace{0.2cm}
\textstyle \int_0^T (Y_s^{i,1} -\displaystyle \max_{j\in \mathcal{I}^{-i}}(Y_s^{j,1} -g_{ij}(s,X_s^{t,x}))) dK_s^{i,1} = 0.
\end{cases}
 \end{split}
\end{equation}
Next, let us set, for $\ii$,
\begin{equation*}
    \underbar Y^i = \bar Y,\quad \underbar Z^i = \bar Z,\quad \underbar V^i = 0,\quad \mbox{ and }\quad \underbar K^i =0.
\end{equation*}
Therefore, $(\underbar Y^i, \underbar Z^i, \underbar V^i, \underbar K^i)_{\ii}$ is a solution of the following system: $\forall \ii$ and $s\leq T$,
\begin{equation*}
\begin{split}
\begin{cases}
\vspace{0.2cm}
\underbar Y^{i}_s = \bar C + \int_s^T \big\lbrace \bar C+ m\,C_f^y\,\underbar Y_r^i + 2\theta \underbar Y_r^{i}+ C_f^z\,|\underbar Z_r^i| \big\rbrace dr\, \\ \vspace{0.2cm}
 \qquad \qquad \quad \,  +\underbar K_T^{i}- \underbar K_s^{i}-\int_s^T \underbar Z_r^{i}dB_r -\int_s^T \int_E \underbar V_r^{i}(e) \tilde{\mu}(dr,de);\\\vspace{0.2cm}
\underbar Y_s^{i} \geqslant  \displaystyle \max_{j \in \mathcal{I}^{-i}}(\underbar Y_s^{j} -g_{ij}(s,X_s^{t,x})); \\\vspace{0.2cm}
\textstyle \int_0^T (\underbar Y_s^{i} -\displaystyle \max_{j\in \mathcal{I}^{-i}}(\underbar Y_s^{j} -g_{ij}(s,X_s^{t,x}))) d\underbar K_s^{i} = 0.
\end{cases}
 \end{split}
\end{equation*}
On the other hand, let $\vec{\Gamma} := (\Gamma^i)_{i=1,...,m} \in \mathcal{H}^{2,m}$ and let us consider the following mapping:
\begin{align}\label{defteta1}
\Theta : \hspace{0.1cm} & \mathcal{H}^{2,m} \rightarrow \mathcal{H}^{2,m} \cr
&\vec{\Gamma} \quad \mapsto \Theta(\vec{\Gamma}):=(Y^{\Gamma,i})_{i=1,...,m}
\end{align}
where $(Y^{\Gamma,i})_{\ii}$ verifies:  $\forall \ii$ and $s\leq T$,
\begin{equation}\label{3.18}
\begin{split}
\begin{cases}
\vspace{0.2cm}
Y^{\Gamma,i}_s = h_i(X_T^{t,x}) +\int_s^T \bar f_i(r,\xtx_r,\vec \Gamma_r,Z^{\Gamma,i}_r, 0)dr +K_T^{\Gamma,i}- K_s^{\Gamma,i}\\ \vspace{0.2cm}
 \qquad \qquad  -\int_s^T Z_r^{\Gamma,i}dB_r -\int_s^T \int_E V_r^{\Gamma,i}(e) \tilde{\mu}(dr,de);\\\vspace{0.2cm}
Y_s^{\Gamma,i} \geqslant  \displaystyle \max_{j \in \mathcal{I}^{-i}}(Y_s^{\Gamma,j} -g_{ij}(s,X_s^{t,x})); \\\vspace{0.2cm}
\textstyle \int_0^T (Y_s^{\Gamma,i} -\displaystyle \max_{j\in \mathcal{I}^{-i}}(Y_s^{\Gamma,j} -g_{ij}(s,X_s^{t,x}))) dK_s^{\Gamma,i} = 0.
\end{cases}
 \end{split}
\end{equation}
 As $\Theta$ is a contraction in $\mathcal{H}^{2,m}$ equipped with  an appropriate equivalent norm (see Proposition 3.3 in \cite{chassagneux2011note}),  then it has a unique fixed point  $(Y^{i,1})_{i\in \mi}$ which, combined with the associated processes $(Z^{i,1},V^{i,1},K^{i,1})_{i\in \mi}$, makes that $(Y^{i,1},Z^{i,1},V^{i,1},K^{i,1})_{i\in \mi}$ is the unique solution of system \eqref{3.14}.\\
 \textbf{}\\
 Now, let us consider the following sequence of processes $(( Y_k^i, Z_k^i,V_k^i, K_k^i)_{\ii})_{k\ge 1}$:
$$
 Y_0^i=0, \mbox{ for all }\ii \mbox{ and for }k\geq 1,\,\,( Y_k^i)_{\ii}=\Theta(( Y_{k-1}^i)_{\ii}),
$$
where $\Theta$ is the mapping defined in \eqref{defteta1}  and $ Z_k^i,\, V_k^i,\, K_k^i$ are associated with $ Y_k^i$, $\ii$, through equation \eqref{3.18}. Therefore, as $\Theta$ is a contraction, the sequence $(( Y_k^i)_{\ii})_{k\ge 0}$ converges to
$(Y^{i,1})_{\ii}$ in $\mathcal{H}^{2,m}$.  On the other hand by an induction argument on $k$ and by using the comparison result, we have that:
\begin{equation}\label{recurrence1}\forall k\geq 0, \,\forall \ii,\,\quad -\bar Y\leq  Y_k^i\leq \bar Y.
\end{equation}
In fact, for $k=0$, this obviously holds since $\bar Y\geq 0$. Next suppose that \eqref{recurrence1} holds for some $k-1$ with $k\geq 1$, i.e.
\begin{equation*}
   \fii,  -\bar Y\leq  Y_{k-1}^i\leq \bar Y.
\end{equation*}
Then, by a linearization procedure of $\bar f_i$,  which is possible since it is Lipschitz $\wt$ $(\vec y,z)$, and  using the induction hypothesis, we obtain: for any $\ii$,
\begin{align*}
    &\bar f_i(s,\xtx_s,(Y_{k-1}^a(s))_{a\in {\cal I}},z,0) \leq \bar C + C_f^y\, \sum_{a=1,m}|Y_{k-1}^a(s)| + C_f^z\,|z|\\
    & \qquad \qquad \qquad\qquad \qquad \quad \quad \leq \bar C +m\,C_f^y\,\bar Y_s + C_f^z\,|z|, \\&\mbox{ and }\\&
    \bar f_i(s,\xtx_s,(Y_{k-1}^a(s))_{a\in {\cal I}},z,0) \geq -(\bar C +m\,C_f^y\,\bar Y_s + C_f^z\,|z|). 
\end{align*}
Finally, by the comparison result (see Proposition 4.2 in \cite{hamadene2015viscosity}) (this is possible since the generators of the systems do not depend on the jump parts), one deduces
that: $\forall \ii$,
\begin{equation*}
     -\bar Y=-\underbar Y^i\leq  Y_{k}^i\leq \underbar Y^i=\bar Y.
\end{equation*}
Taking the limit w.r.t. $k$, we get: $\forall \ii$,
\begin{equation*}
     -\bar Y\leq  Y^{i,1}\leq \bar Y.
\end{equation*}
But, for any $s\in [t,T]\,$ $Y_s^{i,1}=u^{i,1}(s,X_s^{t,x})$. Then, by taking $s=t$ in the previous inequalities, we obtain:
\begin{equation*}
  \forall (t,x)\in [0,T]\times\R^k,\quad |u^{i,1}(t,x)|\leq \bar{Y}_t,
\end{equation*}
which implies that the inequality \eqref{induction} is true for $n=1$. Now, suppose that its holds for some $n-1$ with $n\geq 1$, i.e.,
\begin{equation}\label{3.19}\forall \ii,\,
   \vspace{0.2cm}\forall (t,x)\in [0,T]\times\R^k,\quad |u^{i,n-1}(t,x)|\leq \bar{Y}_t. 
\end{equation}
We are going to prove that, for any $\ii$, $(t,x)\in [0,T]\times \R^k\,$, $|u^{i,n}(t,x)|\leq \bar{Y}_t$.\\
\textbf{}\\
Recall that $(Y^{i,n},Z^{i,n},V^{i,n},K^{i,n})_{i\in \mi}$ the solution of \eqref{recc2} and let us introduce the following mapping:
\begin{align}\label{defteta12}
\Theta^{v,n-1} : \hspace{0.1cm}  \mathcal{H}^{2,m}\quad &\longrightarrow \mathcal{H}^{2,m} \cr
\vec \rho:=(\rho^i)_{\ii}  &\longmapsto \Theta^{v,n-1}(\vec \rho):=(Y^{i,n,\r})_{\ii}
\end{align}
where $(Y^{i,n,\r})_{\ii}$ verifies: $\forall s\le T$,
\begin{equation}\label{recc21}
\begin{split}
\begin{cases}
\vspace{0.2cm}
Y^{i,n,\r}_s = h_i(X_T^{t,x}) +\int_s^T \bar f_i(r,\xtx_r,(\r^k_r)_{\ki},Z^{i,n,\r}_r, \int_E V^{i,n-1}_r(e)\gamma_i(X_r^{t,x},e)\lambda(de))dr\, \\ \vspace{0.2cm}
 \qquad \qquad   +K_T^{i,n,\r}- K_s^{i,n,\r}-\int_s^T Z_r^{i,n,\r}dB_r -\int_s^T \int_E V_r^{i,n,\r}(e) \tilde{\mu}(dr,de);\\\vspace{0.2cm}
Y_s^{i,n,\r} \geq  \displaystyle \max_{j \in \mathcal{I}^{-i}}(Y_s^{j,n,\r} -g_{ij}(s,X_s^{t,x})); \\\vspace{0.2cm}
\textstyle \int_0^T (Y_s^{i,n,\r} -\displaystyle \max_{j\in \mathcal{I}^{-i}}(Y_s^{j,n,\r} -g_{ij}(s,X_s^{t,x}))) dK_s^{i,n,\r} = 0.
\end{cases}
 \end{split}
\end{equation}
Note that $(Y^{i,n})_{\ii}$ verifies $(Y^{i,n})_{\ii}=\Theta^{v,n-1}((Y^{i,n})_{\ii})$ and it is the unique fixed point of 
$\Theta^{v,n-1}$ in $\mathcal{H}^{2,m}$ equipped with  an appropriate equivalent norm (see Proposition 3.3 in \cite{chassagneux2011note}). Next, let us consider the following sequence of processes $((Y^{i,n,l},Z^{i,n,l},V^{i,n,l},K^{i,n,l})_{\ii})_{l\geq 0}$:  
\begin{equation*}
   Y^{i,n,0}=0\mbox{ for all } \ii, \mbox{ and for }\, l\geq 1,\, (Y^{i,n,l})_{\ii}=\Theta^{v,n-1}((Y^{i,n,l-1})_{\ii}),
\end{equation*}
where $(Y^{i,n,l})_{\ii}$ verifies: $\forall s\le T$,
\begin{equation*}
\begin{split}
\begin{cases}
\vspace{0.2cm}
Y^{i,n,l}_s = h_i(X_T^{t,x}) +\int_s^T \bar f_i(r,\xtx_r,(Y^{k,n,l-1}_r)_{k\in \mi},Z^{i,n,l}_r, \int_E V^{i,n-1}_r(e)\gamma_i(X_r^{t,x},e)\lambda(de))dr\, \\ \vspace{0.2cm}
 \qquad \qquad   +K_T^{i,n,l}- K_s^{i,n,l}-\int_s^T Z_r^{i,n,l}dB_r -\int_s^T \int_E V_r^{i,n,l}(e) \tilde{\mu}(dr,de);\\\vspace{0.2cm}
Y_s^{i,n,l} \geqslant  \displaystyle \max_{j \in \mathcal{I}^{-i}}(Y_s^{j,n,l} -g_{ij}(s,X_s^{t,x})); \\\vspace{0.2cm}
\textstyle \int_0^T (Y_s^{i,n,l} -\displaystyle \max_{j\in \mathcal{I}^{-i}}(Y_s^{j,n,l} -g_{ij}(s,X_s^{t,x}))) dK_s^{i,n,l} = 0.
\end{cases}
 \end{split}
\end{equation*}
Similarly as previously, since $\Theta^{v,n-1}$ is a contraction, then the sequence $(( Y^{i,n,l})_{\ii})_{l\ge 0}$ converges to
$(Y^{i,n})_{\ii}$, as $l\rightarrow \infty$, in $\mathcal{H}^{2,m}$. Next, by an induction argument on $l$ and by using the comparison result, we have that:
\begin{equation}\label{3.20}
    \forall \ii,\,\,-\bar Y\leq Y^{i,n,l} \leq \bar Y.
\end{equation}
Actually for $l=0$, the property holds true and if we assume that it is satisfied for some $l-1$ and by using the induction hypotheses, we deduce: $\forall \ii$, $\forall \stt$,
\begin{align*}
    &|\bar f_i(s,\xtx_s,(Y^{k,n,l-1}_s)_{k\in \mi},z, \textstyle \int_E V^{i,n-1}_s(e)\gamma_i(X_s^{t,x},e)\lambda(de))|\\[4pt]
    & = |\bar f_i(s,\xtx_s,(Y^{k,n,l-1}_s)_{k\in \mi},z, \textstyle \int_E \lbrace u^{i,n-1}(s,X_{s_-}^{t,x}+\beta(X_{s_-}^{t,x},e))\\[4pt]
    & \qquad \quad -u^{i,n-1}(s,X_{s_-}^{t,x})\rbrace\gamma_i(X_s^{t,x},e)\lambda(de))|\\[4pt]
    &\leq \bar C + C_f^y\, \displaystyle{\sum_{i=1,m}}|Y^{i,n,l-1}_s| + C_f^z\,|z| +2\theta \bar Y_s\\[4pt]
    & \leq \bar C +m\,C_f^y\,\bar Y_s + C_f^z\,|z| +2\theta \bar Y_s.
\end{align*}
Now by using comparison (Proposition 4.2 in \cite{hamadene2015viscosity}) we deduce that 
\begin{equation*}
    \forall \ii,\,\,-\bar Y=-\underbar Y^i\leq  Y^{i,n,l}\leq \underbar Y^i=\bar Y
\end{equation*}and in taking the limit w.r.t $l\rightarrow\infty$ we obtain: 
\begin{equation*}
    \forall \ii,\,\,-\bar Y\leq Y^{i,n} \leq \bar Y.
\end{equation*}
Finally, as $Y_s^{i,n}=u^{i,n}(s,X_s^{t,x}),\,\, \forall s\in [t,T]$, then 
\begin{equation*}
    \forall (t,x)\in [0,T]\times \R^k,\,\, |u^{i,n}(t,x)|\leq \bar Y_t\leq \underbar C,
\end{equation*}
which implies that $(u^{i,n})_{\ii}$, $n\ge 0$, are uniformly bounded. The proof of the claim is now completed.

Next recall the inequality \eqref{estimationu}. Let us choose $\eta$ a constant such that $ \frac{4}{\alpha_0}Cm\lambda(E)(e^{\alpha_0 \eta }-1)= \frac{3}{4}$. Note that $\eta$ does not depend on the terminal conditions $(h_i)_{\ii}$. Finally let us set
$$\|u^{i,n}- u^{i,p}\|_{\infty,\eta}:=\sup_{(t,x)\in [T-\eta, T]\times \R^k}|u^{i,n}(t,x) - u^{i,p}(t,x)|.
$$
From \eqref{estimationu}, after summation over $i$, we obtain for any $n,p\ge 1$,
\begin{align*}
&\sum_{i=1,m}\|u^{i,n}- u^{i,p}\|^2_{\infty,\eta} \\
& \qquad \leq  \underbrace{\frac{4}{\alpha_0}Cm\lambda(E) (e^{\alpha_0 \eta }-1)}_{=\frac{3}{4}} \sup_{(t,x)\in [T-\eta, T]\times \R^k}\sum_{i=1,m}|u^{i,n-1}(t,x) - u^{i,p-1}(t,x)|^2\\
&\qquad  =  \frac{3}{4}\sum_{i=1,m}\|u^{i,n-1}- u^{i,p-1}\|^2_{\infty,\eta}
\end{align*}
which means that the sequence $((u^{i,n})_{\ii})_{n\geq 0}$ is uniformly convergent in $[T-\eta,T]\times \R^k$. 
Next, let $t\in [T-2\eta, T-\eta]$, then once more by \eqref{estimationu}, we have:
\begin{equation}
\begin{aligned}
&|u^{i,n}(t,x) - u^{i,p}(t,x)|^2 \\[3pt]
& \qquad \leq  2C\E\Big [ \int_{T-2\eta}^{T-\eta}e^{\alpha_0 (r-t)}\int_E \sum_{k=1,m}\{\big| (u^{k, n-1}-u^{k, p-1}) (r,\xtx_{r-}+ \beta(\xtx_{r-},e))\big|^2\\[3pt]
&\quad \qquad + \big|(u^{k,n-1}- u^{k,p-1})(r,\xtx_{r-})\big|^2\} \lambda(de)dr\Big]\\[3pt]
& \quad\qquad  + 2C\E\Big[ \int_{T-\eta}^Te^{\alpha_0 (r-t)}\int_E \sum_{k=1,m}\{\big| (u^{k, n-1}-u^{k, p-1}) (r,\xtx_{r{-}}+ \beta(\xtx_{r-},e))\big|^2\\[3pt] 
& \quad \qquad + \big|(u^{k,n-1}- u^{k,p-1})(r,\xtx_{r{-}})\big|^2\} \lambda(de)dr\Big].
\end{aligned}
\end{equation}
Then, if we set 
$$\|u^{i,n}- u^{i,p}\|_{\infty,2\eta}:=\sup_{(t,x)\in [T-2\eta, T-\eta]\times \R^k}|u^{i,n}(t,x) - u^{i,p}(t,x)|,
$$
we have:
 \begin{align*}
 \sum_{i=1,m}\|u^{i,n}- u^{i,p}\|^2_{\infty,2\eta} &\leq  \frac{4}{\alpha_0}Cm\lambda(E) \Big((e^{\alpha_0 \eta }-1) \sum_{i=1,m}\|u^{i,n-1}- u^{i,p-1}\|^2_{\infty,2\eta}\\[3pt]
 &\qquad  + (e^{2\alpha_0 \eta}-e^{\alpha_0 \eta})\sum_{i=1,m}\|u^{i,n-1}- u^{i,p-1}\|^2_{\infty,\eta}\Big)\\[3pt]
 & \leq  \frac{3}{4}\sum_{i=1,m}\|u^{i,n-1}- u^{i,p-1}\|^2_{\infty,2\eta}\\[3pt]
 &\qquad  +\frac{4}{\alpha_0}Cm\lambda(E)(e^{2\alpha_0 \eta}-e^{\alpha_0 \eta})\sum_{i=1,m}\|u^{i,n-1}- u^{i,p-1}\|^2_{\infty,\eta}.
 \end{align*}
As $\displaystyle \limsup_{n,p\rightarrow \infty}\sum_{i=1,m}\|u^{i,n-1}- u^{i,p-1}\|^2_{\infty,\eta}=0,$ we obtain:
$$
\limsup_{n,p\rightarrow \infty}\sum_{i=1,m}\|u^{i,n}- u^{i,p}\|^2_{\infty,2\eta}\leq  \frac{3}{4}\limsup_{n,p\rightarrow \infty}\sum_{i=1,m}\|u^{i,n-1}- u^{i,p-1}\|^2_{\infty,2\eta}.$$
 Therefore $$
\limsup_{n,p\rightarrow \infty}\sum_{i=1,m}\|u^{i,n}- u^{i,p}\|^2_{\infty,2\eta}=0.
$$
Thus, the sequence $((u^{i,n})_{\ii})_{n\geq 0}$ is uniformly convergent in $[T-2\eta,T-\eta]\times \R^k$. Continuing now this reasoning as many times as necessary on
$[T-3\eta,T-2\eta]$, $[T-4\eta ,T-3\eta]$ and so on, we obtain the uniform convergence of $((u^{i,n})_{\ii})_{n\geq 0}$ in $[0, T]\times \R^k$.
So for  $\ii$ and $(t,x) \in [0,T]\times \R^k$, let us set $u^i(t,x) = \lim_{n \rightarrow \infty} u^{i,n}(t,x)$, $\ii$. Note that $(u^i)_{\ii}$ are continuous and bounded functions on $[0,T]\times \R^k$.

\noindent \tbf{Step 4}: Convergence of $(Y^{i,n},Z^{i,n},V^{i,n},K^{i,n})_{n\geq 0}$

\medskip

\noindent We are now ready to study the convergence of the sequences $(Y^{i,n},Z^{i,n},V^{i,n},K^{i,n})_{n\geq 0}$.

\noindent {\it \underline{Convergence of $(Y^{i,n})_{n\ge 0}$ on $[t,T]$}}: For any $\ii$ and $s\in [t,T]$ let us set $Y^i_s=u^i(s,\xtx_s)$. Next let $n\ge 1$, then:
 \begin{equation}\label{T1}
 \begin{aligned}
 \E\big[\sup_{t\leq s \leq T} |Y_s^{i,n} - Y_s^{i,t,x}|^2 \big]& = \E\big[\sup_{t \leq s \leq T} |{u}^{i,n}(s,X_s^{t,x}) - {u}^{i}(s,X_s^{t,x})\big)|^2 \big] \\[3pt]
 & \leq \|u^{i,n}-u^i\|_\infty:=
 \sup_{(t,x)\in [0, T]\times \R^k}|u^{i,n}(t,x) - u^{i}(t,x)|
 \end{aligned}
 \end{equation}
 As the right hand-side converges to $0$ as $n\rightarrow \infty$, then 
 $(Y^{i,n})_{n\ge 0}$ converges to $(Y^i_s)_{s\in [t,T]}$ in $\ss_{[t,T]}$ which is the space $\ss$ reduced to $[t,T]$. The same is valid for $\aa_{[0,t]}$ (which is $\aa$ reduced to $[0,t]$)(see \eqref{recurrence2} below).
 
 \noindent {\it \usl{Convergence of $(Y^{i,n})_{n\ge 0}$ on $[0,t]$}}: By Remark \ref{surot}, on the time interval $[0,t]$ the sequences
 $(Y^{i,n},Z^{i,n},V^{i,n},K^{i,n})_{n\geq 0}$ verify: 
\begin{equation*}
   \vspace{0.1cm} (Y^{i,0},Z^{i,0},V^{i,0},K^{i,0}) = (0,0,0,0)\, \, \mbox{ for all } \ii, \, \,  \mbox{ for } n\ge 1 \mbox{ and } s \leq t,
  \end{equation*}
\begin{equation}\label{recurrence2}
\begin{split}
\begin{cases}
\vspace{0.2cm} Y^{i,n}\in \ss_{[0,t]}  \mbox{ and } K^{i,n}\in \aa_{[0,t]};\\ \vspace{0.2cm}
Y^{i,n}_s = u^{i,n}(t,x)+\int_s^t \bar f_i(r,x,(Y^{k,n}_r)_{k\in \mi},0,0)dr+K_t^{i,n}- K_s^{i,n},\,\,s\le t;\\\vspace{0.2cm}
Y_s^{i,n} \geqslant  \displaystyle \max_{j \in \mathcal{I}^{-i}}(Y_s^{j,n} -g_{ij}(s,x)),\,\,s\le t; \\\vspace{0.2cm}
\textstyle \int_0^t (Y_s^{i,n} -\displaystyle \max_{j\in \mathcal{I}^{-i}}(Y_s^{j,n} -g_{ij}(s,x))) dK_s^{i,n} = 0.
\end{cases}
 \end{split}
\end{equation}
But $(Y^{i,n}_s)_{s\le t}$ is deterministic, continuous and still have the representation property \eqref{eq9} in connection with the switching problem on $[0,t]$. Next in considering $(\check{P}^{n,p,\d}, \check N^{n,p,\d})$ the solution on $[0,t]$ of the BSDE \eqref{eq8} with generator $\check F^{\d,n,p}(r,x):=
\bar f^\d(r,x,(Y^{k,n}_r)_{k\in \mi},0,0)\vee 
\bar f^\d(r,x,(Y^{k,p}_r)_{k\in \mi},0,0)$
and terminal value $h^{\d,n,p}(t,x):=u^{\d,n}(t,x)\vee u^{\d,p}(t,x)$ and arguing as in Step 3, we deduce a similar inequality as  \eqref{eq9x} that reads: $\forall s\in [0,t]$, 
\begin{align*}\lb{eq9y}
  &\E \Big[  e^{\alpha_0 s}|Y_s^{i,n}- Y_s^{i,p}|^2\Big]\\[4pt]
  &  \leq 2\sum_{k=1,m}e^{\alpha_0 t}|u^{k,n}(t,x)-u^{k,p}(t,x)|^2+ \frac{6C^2}{\alpha_0}\E\Big[ \textstyle \int_s^t  e^{\alpha_0 r}|(Y_r^{k,n})_{k\in \mi}-(Y_r^{k,p})_{k\in \mi}|^2dr \Big].
\end{align*}
Note that the functions $(u^i(t,x))_{\ii}$ verify the consistency condition \eqref{condh} at the terminal time $t$. As we know that for any $\ii$, the sequence 
$(u^{i,n}(t,x))_{n\ge 0}$ converges to $u^i(t,x)$ then it is enough to mimic the arguments of Step 3 to obtain that 
$(Y^{i,n}_s)_{s\le t}$ converges uniformly on $[0,t]$ to some continuous deterministic function (and then bounded) $(Y^{i}_s)_{s\le t}$. 

\noindent As a consequence, for any $\ii$, the sequence $(Y^{i,n})_{n\ge 0}$ converges in $\ss$ to some process $Y^i$, which is moreover bounded since the functions $(u^i)_{\ii}$ are bounded.

\noindent Next as the measure $\lambda$ is finite, then by Remark \ref{surot}, the characterization \eqref{repV} on $[t,T]$ of the sequence $(V^{i,n})_{n \geq 0}$ by means of the function $(u^{i,n})_{n \geq 0}$ and the uniform convergence of $(u^{i,n})_{n \geq 0}$ we deduce that the sequence $(V^{i,n})_{n \geq 0}$ converges in $\hld$ to some process $V^{i,t,x}$ 
which has the following representation:
\begin{align}
V_s^{i}(e) := \{u^{i }(s,\xtx_{s{-}}+& \beta(\xtx_{s{-}},e))- u^{i}(s,\xtx_{s{-}})\}1_{\{s\ge t\}},\,\,ds\otimes d\mathbb{P} \otimes d\lambda \mbox{ on } [0,T] \times \Omega \times E.
\end{align}
This representation imply that $V^{i,t,x}$ are uniformly bounded.
\noindent We now focus on the convergence of the components  $(Z^{i,n}, K^{i,n})_{n\geq 0}$. For this, we first establish a priori estimates, uniform on $n$ of the sequences $(Z^{i,n}, K^{i,n})_{n\geq 0}$. Applying It\^o's formula to $|Y_s^{i,n}|^2$, we have: $\forall s\in [0,T]$
 \begin{equation}
\begin{aligned}
\E \big[|Y_s^{i,n}|^2 \big] &+ \E \big[\int_s^T|Z_r^{i,n}|^2dr \big] + \E \big[\int_s^T\int_E|V_r^{i,n}(e)|^2 \lambda(de)dr \big] \\[6pt] 
&= \E \big[|h_i(X_T^{t,x}|^2 \big] + 2\E \big[ \int_s^T Y_r^{i,n} f_i(r,X_r^{t,x},(Y_r^{k,n})_{k\in \mi}, Z_r^{i,n}, V_{r}^{i,n-1})dr\big]
+ 2 \E \big[ \int_s^T Y_r^{i,n} dK_r^{i,n}\big].
\end{aligned}
\end{equation}
Then by a linearization procedure of $f_i(.)$,  which is possible since it is Lipschitz $\wt$ $(\vec y,z,q)$ and using the inequality $2ab\leq \frac{1}{\epsilon} a^2 + \epsilon b^2$ for any constant $\epsilon > 0$, we have:
 \begin{equation*}
\begin{aligned}
 &\E \big[\textstyle \int_0^T|Z_r^{i,n}|^2dr \big]\\[4pt]
 & \quad \leq \E \big[|h_i(X_T^{t,x}|^2 \big]  + 2\E \big[ \textstyle \int_0^T |Y_r^{i,n}|\{|f_i(r,X_r^{t,x},0,0,0)|+ \sum_{l=1,m}a_r^{i,l,n}|Y_r^{l,n}|+  b_r^{i,n}|Z_r^i|\\[4pt] 
& \qquad + c_r^{i,n}\textstyle \int_E |V_{r}^{i,n-1}(e)\g(\xtx_r,e)| \lambda(de)\} dr\big] + \frac{1}{\epsilon}\E \big[ \sup_{s\leq T} |Y_s^{i,n}|^2 \big] + \epsilon \E \big[ \big(K_T^{i,n}\big)^2\big],
\end{aligned}
\end{equation*}

\noindent where $a^{i,l,n}\in \R$, $b^{i,l,n}\in \R^d$  are ${\mathcal P}$-measurable non-negative bounded processes while $c^{i,n}\in \R$ is non-negative bounded and $\mathbf{P}$-measurable. Using again the inequality $2ab\leq \frac{1}{\nu} a^2 + \nu b^2$ for $\nu > 0$, yields 
\begin{equation*}
\begin{aligned}
 &\E \big[\textstyle \int_0^T|Z_r^{i,n}|^2dr \big]\\[4pt]
 & \quad \leq \E \big[|h_i(X_T^{t,x})|^2 \big]  + \frac{1}{\nu}\E \big[ \textstyle \int_0^T |Y_r^{i,n}|^2dr \big]+ \nu  \E \big[ \int_0^T\{ |f_i(r,X_r^{t,x},0,0,0)| \\[4pt] 
& \qquad  +\textstyle \sum_{l=1,m}a_r^{i,l,n}|Y_r^{l,n}|+  b_r^{i,n}|Z_r^i|+ c_r^{i,n}\textstyle \int_E |V_{r}^{i,n-1}(e)\g(\xtx_r,e)| \lambda(de)\}^2 dr\big]\\[4pt]
& \qquad + \frac{1}{\epsilon}\E \big[ \sup_{s\leq T} |Y_s^{i,n}|^2 \big] + \epsilon \E \big[ \big(K_T^{i,n}\big)^2\big].
\end{aligned}
\end{equation*}
From the boundedness of $f_i(t,x,0,0,0)$ and $h_i(x)$, the inequality $|a + b + c + d|^2 \leq 4 \{|a|^2 + |b|^2 + |c|^2 + |d|^2\},$ $\forall a, b, c, d\in \R$ and finally the Cauchy-Schwarz one, we have:
\begin{equation*}
\begin{aligned}
&\E \big[\textstyle \int_0^T|Z_r^{i,n}|^2dr \big] \\[4pt]
& \quad \textstyle \leq \bar C^2+4\nu \bar C^2 T + \frac{1}{\nu}\E \big[ \int_0^T |Y_r^{i,n}|^2dr \big]+  4\nu C_3\E\big[\int_0^T\sum_{l=1,m}|Y_r^{l,n}|^2 dr \big]\\[4pt]
 & \qquad +4\nu C_3\E\big[\textstyle \int_0^T |Z_r^{i,n}|^2dr \big] + 4\nu C_3\E\big[\int_0^T\int_E |V_{r}^{i,n-1}(e)|^2 \lambda(de) dr\big]\\[6pt] 
& \qquad + \textstyle \frac{1}{\epsilon}\E \big[ \sup_{s\leq T} |Y_s^{i,n}|^2 \big] + \epsilon \E \big[ \big(K_T^{i,n}\big)^2\big],
\end{aligned}
\end{equation*}
 for a suitable positive constants $C_1, C_2$ and $C_3$. Choose now $\nu$ such that $4\nu C_3 < 1$, and taking the sum over all $i \in \mi$, we obtain:
 \begin{equation*}
\begin{aligned}
 \sum_{i=1,m}&\E \big[\int_0^T|Z_r^{i,n}|^2dr \big]\\[4pt]
 & \quad \leq C\Big( 1 + \textstyle \sum_{i=1,m}\E \big[\sup_{s \leq T}|Y_s^{i,n}|^2 \big]  + \sum_{i=1,m}\E\big[\int_0^T\int_E |V_{r}^{i,n-1}(e)|^2 \lambda(de) dr\big]\Big)+ \epsilon\textstyle \sum_{i=1,m}\E \big[ \big(K_T^{i,n}\big)^2\big],
\end{aligned}
\end{equation*} 
where $C = C(T,m,\nu,\epsilon) >0$ is an appropriate constant independent of $n$. Thanks to the convergence of $(Y^{i,n})_n$  in $\ss$, we have $\sup_{n \geq 0} \E \big[\sup_{s \leq T}|Y_s^{i,n}|^2\big] \leq C$, and then taking into consideration the convergence of  $(V^{i,n})_n$  in $\hld$, we finally obtain
\begin{equation}\label{estimZ}
\begin{aligned}
 \textstyle \sum_{i=1,m}\E \big[\int_0^T|Z_r^{i,n}|^2dr \big] & \leq C+ \epsilon \textstyle \sum_{i=1,m}\E \big[ \big(K_T^{i,n}\big)^2\big].
\end{aligned}
\end{equation} 
Now, from the relation 
\begin{equation}\label{relationK}
\begin{aligned}
K_T^{i,n} = &\, Y_0^{i,n} - h_i(X_T^{t,x}) - \textstyle\int_0^T f_i(r,X_r^{t,x},(Y_r^{k,n})_{k\in \mi}, Z_r^{i,n}, V_{r}^{i,n-1})dr \\[4pt]
& \qquad +\textstyle \int_0^T Z_r^{i,n} dB_r + \int_0^T \int_E V_r^{i,n}(e)\tilde{\mu}(dr,de),
\end{aligned}
\end{equation}
and once again, by the linearization procedure of the Lipschitz function $\bar f_i(.)$ and the boundedness of $f_i(t,x,0,0,0)$ and $h_i(x)$, there exist some positive constant $C^{\prime}$ such that
 \begin{equation*}
\begin{aligned}
 \textstyle \sum_{i=1,m}\E \big[\big(K_T^{i,n}\big)^2 \big] & \leq C^{\prime} \Big( 1 +\textstyle \sum_{i=1,m}\E \big[\sup_{s \leq T}|Y_s^{i,n}|^2 \big] + \sum_{i=1,m}\E \big[\int_0^T|Z_r^{i,n}|^2dr \big]\\[4pt]
 & \qquad \qquad + \textstyle \sum_{i=1,m}\E\big[\int_0^T\int_E |V_{r}^{i,n-1}(e)|^2 \lambda(de) dr\big] \Big)\\[4pt]
 & \qquad \leq C^{\prime} \Big( 1 + \textstyle \sum_{i=1,m}\E \big[\int_0^T|Z_r^{i,n}|^2dr \big] \Big).
\end{aligned}
\end{equation*} 
Combining this last estimate with \eqref{estimZ} and choosing $\epsilon$ small enough since it is arbitrary, there exists a constant $\bar{C}$ such that 
 \begin{equation}\label{estimZK}
\begin{aligned}
 \textstyle \sum_{i=1,m}\E \big[\int_0^T|Z_r^{i,n}|^2dr + \big(K_T^{i,n}\big)^2  \big] & \leq \bar{C}.
\end{aligned}
\end{equation}
Next, for any $n, p \geq 1$, by It\^o's formula we have:
\begin{equation*}
\begin{aligned}
 &\E \big[\textstyle \int_0^T|Z_r^{i,n} -Z_r^{i,p} |^2dr \big]\\[4pt] 
& \quad \leq 2\E \big[\textstyle \int_0^T \big(Y_r^{i,n} -Y_r^{i,p}\big) \big(f_i(r,X_r^{t,x},(Y_r^{k,n})_{k\in \mi}, Z_r^{i,n}, V_{r}^{i,n-1})- f_i(r,X_r^{t,x},(Y_r^{k,p})_{k\in \mi}, Z_r^{i,p}, V_{r}^{i,p-1})\big)dr\big]\\[4pt]
& \qquad\qquad +2 \E \big[ \textstyle\int_0^T \big(Y_r^{i,n} -Y_r^{i,p}\big)  \big(dK_r^{i,n}(r) -dK_r^{i,p}\big) \big].
\end{aligned}
\end{equation*}
By Cauchy-Schwarz inequality and using the inequality $2ab\leq \frac{1}{\eta} a^2 + \eta b^2$ for $\eta > 0$, we have:
\begin{equation*}
\begin{aligned}\textstyle 
 &\E \big[\textstyle \int_0^T|Z_r^{i,n} -Z_r^{i,p}|^2dr \big] 
\leq  2 \, \sqrt{\E \big[ \sup_{s \leq T} |Y_r^{i,n} -Y_r^{i,p}|^2\big]} \times \\[3pt]
 & \sqrt{\E \big[\textstyle  \int_0^T \big|f_i(r,X_r^{t,x},(Y_r^{k,n})_{k\in \mi}, Z_r^{i,n}, V_{r}^{i,n-1})
- f_i(r,X_r^{t,x},(Y_r^{k,p})_{k\in \mi}, Z_r^{i,p}, V_{r}^{i,p-1})\big|^2dr\big]}\\[3pt]
&  \qquad +  \frac{1}{\eta}\,\E \big[ \sup_{s\leq T} |Y_s^{i,n} - Y_s^{i,p}|^2 \big] + \eta \, \E \big[ \big(K_T^{i,n}+ K_T^{i,p}\big)^2\big].
\end{aligned}
\end{equation*}
But there exists a constant $C \geq 0$ (independent of $n$ and $p$) such that, for all $n,p \geq 1$,
\begin{equation}\label{estimationf}
\begin{aligned}
 \E \big[ \textstyle \int_0^T \big|f_i(r,X_r^{t,x},(Y_r^{k,n})_{k\in \mi}, Z_r^{i,n}, V_{r}^{i,n-1})
- f_i(r,X_r^{t,x},(Y_r^{k,p})_{k\in \mi}, Z_r^{i,p}, V_{r}^{i,p-1})\big|^2dr\big] \leq C. 
\end{aligned}
\end{equation}
Then taking the limit w.r.t $n,p$ in the previous inequality and taking into account of \eqref{estimZK} and the convergence of  $Y^{i,n}$ in $\ss$, we deduce that:
\begin{equation*}
\begin{aligned}
 \limsup_{n,p\rightarrow \infty}\E &\big[\textstyle \int_0^T|Z_r^{i,n} -Z_r^{i,p} |^2dr \big] \le \bar C \eta. \end{aligned}
\end{equation*}
As $\eta$ is arbitrary then $(Z^{i,n})_{n \geq 0}$ is a Cauchy sequence in $\hdd$. Therefore there exists a process $Z^{i,t,x}$ of $\hdd$ such that $(Z^{i,n})_{n \geq 0}$ converges to $Z^{i,t,x}$ in $\hdd$. Finally, since for $s\leq T,$
\begin{equation*}
\begin{aligned}
 K_s^{i,n} = &\, Y_0^{i,n} -Y_s^{i,n}- \textstyle \int_0^s f_i(r,X_r^{t,x},(Y_r^{k,n})_{k\in \mi}, Z_r^{i,n}, V_r^{i,n-1})dr \\[3pt]
& \qquad +\textstyle  \int_0^s Z_r^{i,n} dB_r + \int_0^s \int_E V_r^{i,n}(e)\tilde{\mu}(dr,de),  
\end{aligned}
\end{equation*}
then, we have also $ \E \big[\sup_{s \leq T}|K_s^{i,n} -K_s^{i,p} |^2 \big] \rightarrow 0 \, \, \mbox{ as } \, \,  n, p\rightarrow \infty$. Thus, there exists a process $(K_s^{i,t,x})_{s \leq T}$ which belongs to $\aa$ such that $\E \big[\sup_{s \leq T}|K_s^{i,n} -K^{i,t,x}_s |^2 \big] \rightarrow 0 \, \, \mbox { as } \, \,  n \rightarrow \infty$. Moreover we have: $\forall \sot$,
\begin{equation}\label{eqRBSDExx}
\begin{split}
\begin{cases}
\vspace{0.3cm}
Y_s^{i,t,x} = h_i(X_T^{t,x})+ \int_s^T \bar f_i(r,X_r^{t,x},(Y_r^{k,t,x})_{k\in \mi},Z_r^{i,t,x}, \int_E V_r^{i,t,x}(e)\gamma_i(X_r^{t,x},e) \lambda(de))dr \\\vspace{0.3cm}
 \qquad \qquad   +K_T^{i,t,x} - K_s^{i,t,x} -\int_s^T Z_r^{i,t,x}dB_r -\int_s^T \int_E V_r^{i,t,x}(e) \tilde{\mu}(dr,de);\\\vspace{0.3cm}
 Y_s^{i,t,x} \geqslant  \displaystyle \max_{j \in \mathcal{I}^{-i}}(Y_s^{j,t,x} -g_{ij}(s,X_s^{t,x})).
\end{cases}
\end{split}
\end{equation} 
Finally, let us show that the third condition in \eqref{eqBSDE} is satisfied by $(Y^{i,t,x},Z^{i,t,x},V^{i,t,x},K^{i,t,x})_{\ii}$. Actually 
\begin{equation}\label{thirdcondition}
\begin{aligned} 
\textstyle \int_0^T (Y_s^{i,t,x} -&\displaystyle \max_{j\in \mathcal{I}^{-i}}(Y_s^{j,t,x} -g_{ij}(s,X_s^{t,x}))) dK_s^{i,t,x}\\[2pt]
&= \textstyle \int_0^T (Y_s^{i,t,x} -\displaystyle \max_{j\in \mathcal{I}^{-i}}(Y_s^{j,t,x} -g_{ij}(s,X_s^{t,x}))) (dK_s^{i,t,x} - dK_s^{i,n})\\[2pt]
& +\textstyle  \int_0^T (Y_s^{i,t,x} -\displaystyle \max_{j\in \mathcal{I}^{-i}}(Y_s^{j,t,x} -g_{ij}(s,X_s^{t,x}))) dK_s^{i,n}. 
\end{aligned}
\end{equation}
Let $\omega$ be fixed. It follows from the uniform convergence of $(Y^{i,n})_n$ to $(Y^{i,t,x})_{\ii}$ that, for any $\epsilon \geq 0$, there exist $N_{\epsilon}( \omega) \in \N$, such that for any $n \geq N_{\epsilon}( \omega)$ and for any $s \leq T$,
\begin{align*}
     Y_s^{i,t,x}(\omega)  -\displaystyle \max_{j\in \mathcal{I}^{-i}}(Y_s^{j,t,x}(\omega) &-g_{ij}(s,X_s^{t,x}(\omega)))\\[2pt]
     & \quad \leq Y_s^{i,n}(\omega) -\displaystyle \max_{j\in \mathcal{I}^{-i}}(Y_s^{j,n}(\omega) -g_{ij}(s,X_s^{t,x}(\omega))) + \epsilon.
\end{align*}
Therefore, for  $n \geq N_{\epsilon}( \omega)$ we have
\begin{equation}\label{thirdcondition2}
 \int_0^T (Y_s^{i,t,x} -\displaystyle \max_{j\in \mathcal{I}^{-i}}(Y_s^{j,t,x} -g_{ij}(s,X_s^{t,x}))) dK_s^{i,n} \leq \epsilon K^{i,n}_T(\omega).
\end{equation}
On the other hand, the function 
\begin{align*}
    Y^{i,t,x}(\omega)- \displaystyle \max_{j\in \mathcal{I}^{-i}}(Y^{j,t,x}(\omega) -g_{ij}(.,\xtx_.(\omega))) :
 s \in [0,T] \longmapsto Y^{i,t,x}_s(\omega)- \displaystyle \max_{j\in \mathcal{I}^{-i}}(Y^{j,t,x}_s(\omega) -g_{ij}(s,X_s^{t,x}(\omega)))
\end{align*}
is c\`adl\`ag and then bounded, then there exists a sequence of step functions $(f^m(\omega,.))_{m\ge 1}$ which converges uniformly on $[0,T]$ to $Y^{i,t,x}(\omega)- \displaystyle \max_{j\in \mathcal{I}^{-i}}(Y^{j,t,x}(\omega) -g_{ij}(.,\xtx_.(\omega)))$, i.e., there exist $m_{\epsilon}(\omega) \geq 0$ such that for $m \geq m_{\epsilon}(\omega)$, we have 
$$\forall s\le T,\,\,\big| Y^{i,t,x}_s(\omega)- \displaystyle \max_{j\in \mathcal{I}^{-i}}(Y^{j,t,x}_s(\omega) -g_{ij}(s,\xtx_s(\omega))) - f^m(\omega,s) \big| < \epsilon.$$ 
It follows that 
\begin{equation*}
\begin{aligned}
&\textstyle \int_0^T (Y_s^{i,t,x} -\displaystyle \max_{j\in \mathcal{I}^{-i}}(Y_s^{j,t,x} -g_{ij}(s,X_s^{t,x}))) (dK_s^{i,t,x} - dK_s^{i,n})\\[3pt]
& =\textstyle \int_0^T (Y_s^{i,t,x} -\displaystyle \max_{j\in \mathcal{I}^{-i}}(Y_s^{j,t,x} -g_{ij}(s,X_s^{t,x})) - f^m(\omega,s)) (dK_s^{i,t,x} - dK_s^{i,n})\\[3pt]
& \quad + \textstyle  \int_0^T f^m(\omega,s)(dK_s^{i,t,x} - dK_s^{i,n})\\[3pt]
& \leq \textstyle\int_0^T  f^m(\omega,s)(dK_s^{i,t,x} - dK_s^{i,n}) + \epsilon (K_T^{i,t,x}(\omega) + K_T^{i,n}(\omega)).
\end{aligned}
\end{equation*}
But the right-hand side converges to $2\epsilon K_T^{i,t,x}(\omega)$, as $n \rightarrow \infty$, since  $f^m(\omega,.)$ is a step function and then\\ $\int_0^T  f^m(\omega,s)(dK_s^{i,t,x} - dK_s^{i,n}) \rightarrow 0$ as $n\rw \infty.$ Therefore, we have 
\begin{equation} \label{thirdconction3}\limsup_{n\rightarrow \infty}\int_0^T (Y_s^{i,t,x} -\displaystyle \max_{j\in \mathcal{I}^{-i}}(Y_s^{j,t,x} -g_{ij}(s,X_s^{t,x}))) (dK_s^{i,t,x} - dK_s^{i,n})\leq 2 \epsilon K_T^{i,t,x}.
\end{equation}
Finally, from \eqref{thirdcondition}, \eqref{thirdcondition2} and  \eqref{thirdconction3} we deduce that 
\begin{equation*} 
\int_0^T (Y_s^{i,t,x}(\omega) -\displaystyle \max_{j\in \mathcal{I}^{-i}}(Y_s^{j,t,x}(\omega) -g_{ij}(s,X_s^{t,x})(\omega))) dK_s^{i,t,x}(\omega) \leq 3 \epsilon K_T^{i,t,x}(\omega).
\end{equation*}
As $\epsilon$ is arbitrary and $Y_s^{i,t,x} \geq \displaystyle \max_{j\in \mathcal{I}^{-i}}(Y_s^{j,t,x} -g_{ij}(s,X_s^{t,x}))$, then
$$ \int_0^T (Y_s^{i,t,x} -\displaystyle \max_{j\in \mathcal{I}^{-i}}(Y_s^{j,t,x} -g_{ij}(s,X_s^{t,x}))) dK_s^{i,t,x} =0,$$
which completes the proof.
\qed
Now, we study the system \eqref{eqBSDE} in the general case, i.e., without assuming the boundedness of the functions 

The following is the main result of this section.  
\begin{theoreme}
Assume that the functions $(\bar f_i)_{i\in \mathcal{I}}$ and $(\gamma_i)_{\ii}$ verify Assumption (H1) and, $(g_{ij})_{i,j \in \mathcal{I}}$  and $(h_i)_{i \in \mathcal{I}}$ verify Assumptions (H2) and (H3).
Then the system \eqref{eqBSDE} has a solution $(Y^{i,t,x},Z^{i,t,x},V^{i,t,x},K^{i,t,x})_{i\in \mi}$.
Moreover there exists continuous functions $(u^i)_{\ii}$ of polynomial growth such that for any $\ii$, $\tx$,
$$Y^{i,t,x}_s=u^i(s,\xtx_s), \,\,\forall s\in [t,T].
$$
\end{theoreme}
\nd \sol{Proof}: First we are going to transform the system \eqref{eqBSDE} in such a way to fall in the same framework as the one of Theorem \ref{thmexistence2}. So let $\phi$ be a function defined as follows ($p$ is the same or greater than the exponents which are involved in (H1)-iii) and (H3)): \begin{equation} \label{phi}
\phi(x):= \frac{1}{(1+|x|^2)^p},\,\,x\in \R^k, 
\end{equation} and for $\sot$ let us define,
\begin{equation}\label{ybar}
\overline{Y}_s^{i,t,x} := Y_s^{i,t,x} \phi(X_s^{t,x}).
\end{equation}
Then, by It\^o's formula we have: $\forall s \in [0,T]$,
\begin{align*}
   \phi(X_s^{t,x})= &\phi(X_0^{t,x}) + \textstyle \int_0^s D_x\phi(X_{r-}^{t,x})dX_r^{t,x} +  \frac{1}{2}\textstyle \int_0^s \mbox{Tr}(D^2_{xx}\phi(X_{r-}^{t,x})\sigma\sigma^\top(r,X_r^{t,x}))dr\\[7pt]
   &+\textstyle  \sum_{0<r\leq s}\{\phi(X_{r}^{t,x}) -\phi(X_{r-}^{t,x}) - D_x\phi(X_{r-}^{t,x})\Delta_r X^{t,x}\}.
\end{align*}
Since $X^{t,x}$ satisfies the SDE \eqref{eq1}, then for $s\in [0,T]$,
\begin{align*}
   \textstyle \sum_{0<r\leq s}\{\phi(X_{r}^{t,x})&-\phi(X_{r-}^{t,x})-D_x\phi(X_{r-}^{t,x})\Delta_r X^{t,x}\}\\[7pt]
    &=\textstyle \sum_{0<r\leq s}\{\phi(X_{r-}^{t,x} + \Delta_rX^{t,x}) -\phi(X_{r-}^{t,x}) - D_x\phi(X_{r-}^{t,x})\Delta_r X^{t,x}\}\\[7pt]
    & = \textstyle \int_0^s\int_E \{\phi(X_{r-}^{t,x} + \beta(X_{r-}^{t,x},e)) -\phi(X_{r-}^{t,x}) - D_x\phi(X_{r-}^{t,x}) \beta(X_{r-}^{t,x},e)\} \mu(dr,de)\\[7pt]
    & = \textstyle \int_0^s\int_E \{\phi(X_{r-}^{t,x} + \beta(X_{r-}^{t,x},e)) -\phi(X_{r-}^{t,x}) - D_x\phi(X_{r-}^{t,x}) \beta(X_{r-}^{t,x},e)\} \tilde \mu(dr,de)\\[7pt]
    & \qquad +\textstyle \int_0^s\int_E \{\phi(X_{r-}^{t,x} + \beta(X_{r-}^{t,x},e)) -\phi(X_{r-}^{t,x}) - D_x\phi(X_{r-}^{t,x}) \beta(X_{r-}^{t,x},e)\} \lambda(de)ds.
\end{align*}
 Next, going back to \eqref{ybar} and using It\^o's formula we obtain: $\forall s \in [0,T]$,
 \begin{equation*}
\begin{aligned}
d\overline{Y}_s^{i,t,x} & =  Y_{s-}^{i,t,x} d\phi(X_s^{t,x}) + \phi(X_{s-}^{t,x})dY_s^{i,t,x} + d[ Y^{i,t,x},\phi(X^{t,x})]_s,
\end{aligned}
\end{equation*}
where 
\begin{align*}
    [ Y^{i,t,x},\phi(X^{t,x})]_s = \langle  Y^{i,t,x}, \phi(X^{t,x})\rangle_s^c + \sum_{0<r\leq s} \Delta_r Y^{i,t,x} \Delta_r \phi(X^{t,x}).
\end{align*}
But
\begin{align*}
    d\langle  Y^{i,t,x}, \phi(X^{t,x})\rangle_s^c =
    Z_s^{i,t,x} D_x\phi(X_{s-}^{t,x}) \sigma(s,X_s^{t,x})ds
\end{align*}
and 
\begin{align*}
    \sum_{0<r\leq s} \Delta_r Y^{i,t,x} \Delta_r \phi(X^{t,x})&= \sum_{0<r\leq s} \Delta_r Y^{i,t,x}\{\phi(X_{r}^{t,x}) -\phi(X_{r-}^{t,x})\} \\[3pt]
    & = \sum_{0<r\leq s} \Delta_r Y^{i,t,x}\{\phi(X_{r-}^{t,x}+ \Delta_rX^{t,x}) -\phi(X_{r-}^{t,x})\} \\[3pt]
    & = \textstyle \int_0^s \int_E V_r^{i,t,x}(e)\{\phi(X_{r-}^{t,x}+ \beta(\xtx_{r-},e) ) -\phi(X_{r-}^{t,x})\} \mu(dr,de)\\[3pt]
     & = \textstyle \int_0^s \int_E V_r^{i,t,x}(e)\{\phi(X_{r-}^{t,x}+ \beta(\xtx_{r-},e) ) -\phi(X_{r-}^{t,x})\} \tilde \mu(dr,de)\\[3pt]
      & \qquad + \textstyle \int_0^s \int_E V_r^{i,t,x}(e)\{\phi(X_{r-}^{t,x}+ \beta(\xtx_{r-},e) ) -\phi(X_{r-}^{t,x})\} \lambda(de)dr.
\end{align*}
Then it follows that: $\forall s \in [0,T]$,
\begin{align*}
  d\overline{Y}_s^{i,t,x}  =  \Big\lbrace -\phi(X_{s-}^{t,x}) &\bar f_i(s,X_s^{t,x},(Y_s^{k,t,x})_{k\in \mi},Z_s^{i,t,x},\textstyle\int_E V^{i,t,x}_s(e)\gamma_i(X_s^{t,x},e)\lambda(de)) \\[7pt]
  & \quad  + Y_{s-}^{t,x}\big(D_x\phi(X_{s-}^{t,x})b(s,X_s^{t,x}) + \frac{1}{2}\tr(D^2_{xx}\phi(X_{s-}^{t,x})  \sigma\sigma^\top(s,X_s^{t,x}))\big)\\[7pt]
  & \quad +  Y_{s-}^{t,x}\textstyle\int_E\big( \phi(X_{s-}^{t,x} +\beta(X_{s-}^{t,x},e)) - \phi(X_{s-}^{t,x}) -D_x\phi(X_{s-}^{t,x})\beta(X_{s-}^{t,x},e)\big)\lambda(de)\\[7pt]
  & \quad  +Z_s^{i,t,x} D_x\phi(X_{s-}^{t,x}) \sigma(s,X_s^{t,x}) +\textstyle \int_E V_s^{i,t,x}(e)\{\phi(X_{s-}^{t,x}+ \beta(\xtx_{s-},e) ) \\[7pt]
   &  \quad -\phi(X_{s-}^{t,x})\}\lambda(de) \Big\rbrace ds- \phi(X_{s-}^{t,x}) dK_s^{i,t,x} + \Big\lbrace \phi(X_{s-}^{t,x}) Z_s^{i,t,x}   \\[7pt]
    & \quad+ Y_{s-}^{i,t,x} D_x\phi(X_{s-}^{t,x})\sigma(s,X_s^{t,x})\Big\rbrace dB_s +\textstyle \int_E \Big\lbrace Y_{s-}^{i,t,x} \{ \phi(X_{s-}^{t,x} + \beta(X_{s-}^{t,x},e))\\[7pt]
    & \quad -\phi(X_{s-}^{t,x})\} +  V_s^{i,t,x}(e)\phi(X_{s-}^{t,x}+ \beta(\xtx_{s-},e))\Big\rbrace \tilde\mu(ds,de).
\end{align*}
Next let us set, for  $s\in [0, T]$,
\begin{align*}
 &d\overline{K}_s^{i,t,x} :=\phi(X_{s-}^{t,x})dK_s^{i,t,x} \mbox{ and }\overline{K}_0^{i,t,x}=0,\\
&\overline{Z}_s^{i,t,x} :=\phi(X_{s-}^{t,x}) Z_s^{i,t,x} + Y_{s-}^{i,t,x} D_x\phi(X_{s-}^{t,x})\sigma(s,X_s^{t,x}) \,\,\, \,  \mbox{ and }\\
&\overline{V}_s^{i,t,x}(e) :=  Y_{s-}^{i,t,x} \{ \phi(X_{s-}^{t,x} + \beta(X_{s-}^{t,x},e))-\phi(X_{s-}^{t,x})\} +  V_s^{i,t,x}(e)\phi(X_{s-}^{t,x}+ \beta(\xtx_{s-},e)).
\end{align*}
Then  $((\overline{Y}^{i,t,x} ,\overline{Z}^{i,t,x} ,\overline{V}^{i,t,x},\overline{K}^{i,t,x} ))_{\ii}$ verifies: $\forall s\in [0,T]$,
\begin{equation}\label{nv}
\begin{cases}
\vspace{0.3cm}\overline{Y}_s^{i,t,x} = \check{h}_i(X_T^{t,x})+\int_s^T\check{F}_i(r,X_r^{t,x},(\overline{Y}_r^{k,t,x})_{k\in \mi},\overline{Z}_r^{i,t,x},\int_E\overline{V}_{r}^{i,t,x}\gamma_i(X_r^{t,x},e) )dr  \\\vspace{0.3cm}
\qquad  \qquad +\overline{K}_T^{i,t,x}- \overline{K}_s^{i,t,x}  -\int_s^T \overline{Z}_r^{i,t,x}dB_r- \int_s^T \int_E \overline{V}_r^{i,t,x}(e)\tilde{\mu}(dr,de),\\\vspace{0.3cm}
\overline{Y}_s^{i,t,x} \geqslant  \displaystyle \max_{j \in \mathcal{I}^{-i}}(\overline{Y}_s^{j,t,x} -\check{g}_{ij}(s,X_s^{t,x}))\\\vspace{0.3cm}
\textstyle \int_0^T (\overline{Y}_s^{i,t,x} -\displaystyle \max_{j\in \mathcal{I}^{-i}}(\overline{Y}_s^{j,t,x} -\check{ g}_{ij}(s,X_s^{t,x}))) d\overline{K}_s^{i,t,x} = 0,
\end{cases}
\end{equation}
where for any $i,j\in \mi$,
$$
\check{h}_i(X_T^{t,x}):= \phi(X_T^{t,x}) h_i(X_T^{t,x}),\,\, \check{g}_{ij}(s,X_s^{t,x}) :=\phi(X_s^{t,x}) g_{ij}(s,X_s^{t,x}), $$
and for any $(s,x,\vy,z,v)\in [0,T]\times \R^{k+m+d}\times L^2(d\lambda)$, 
\begin{align*}
\check{F}_i &(s,x,\vec{y}, z,v) :=\phi(x)\bar f_i \Big [s,x,\phi(x)^{-1}\vec y, \phi(x)^{-1}z-y^i\phi(x)^{-2}D_x\phi(x) \sigma(s,x),\\[5pt]
&\quad \textstyle \int_E \gamma_i(x,e) \phi(x+\beta(x,e))^{-1}v(e)\lambda(de)-y^i
\int_E\gamma_i(x,e) \phi(x+\beta(x,e))^{-1}  \\[5pt]
& \quad \times\phi(x)^{-1}(\phi(x+\beta(x,e))-\phi(x))\lambda(de)
\Big ]-y^i\phi(x)^{-1} \Big \{ b(s,x) D_x\phi(x) \\[5pt]
&  \quad +\frac{1}{2} \tr (D_{xx}^2\phi(x)\sigma \sigma^\top(s,x) )+ \textstyle \int_E (\phi(x+\beta(x,e))-\phi(x) - D_x\phi(x)\beta(x,e))\lambda(de)\\[5pt]
& \quad \textstyle +\int_E\lambda(de)
(\phi(x+\beta(x,e))-\phi(x))^2\phi(x+\beta(x,e))^{-1}
\Big \}\\[5pt]
&\quad  - z\phi(x)^{-1}D_x\phi(x)^\top\sigma(s,x) -y^i\phi(x)^{-2}D_x\phi(x)^\top \sigma(s,x)\sigma(s,x)^\top
D_x\phi(x)\\[5pt]
&\quad -\textstyle \int_E(\phi(x+\beta(x,e))-\phi(x))\phi(x+\beta(x,e))^{-1}
v(e)\lambda(de).
\end{align*}
Here, let us notice that the functions $(\check{g}_{ij})_{i,j\in \mathcal{I}}$ and $(\check{h}_i)_{\ii}$ verify Assumptions (H2)-(H3) while $(\check{F_i})_{\ii}$ satisfy 
(H1)-ii), iii), iv). 

By Theorem \ref{thmexistence1}, the following scheme $(\check{Y}^{i,n,t,x} ,\check{Z}^{i,n,t,x} ,\check {V}^{i,n,t,x},\check{K}^{i,n,t,x} )_{\ii}$, $n\ge 1$, 
is well-posed:
$\check{V}^{i,0,t,x}=0$ and for $n\ge 1$ (we omit the dependence on $t,x$ as there is no confusion) and $\sot$,
\begin{equation}\label{recurrencex2}
\begin{split}
\begin{cases}
\vspace{0.2cm} \check Y^{i,n}\in \ss,  \check Z^{i,n} \in \hdd, \check V^{i,n} \in \hld, \mbox{ and } \check K^{i,n}\in \aa;\\ \vspace{0.2cm}
\check Y^{i,n}_s = \check h_i(X_T^{t,x}) +\int_s^T \check  F_i(r,\xtx_r,(\check Y^{k,n}_r)_{k\in \mi},\check Z^{i,n}_r,\check V^{i,n-1}_r)dr \\ \vspace{0.4cm}
 \qquad \qquad \quad \,  +\check K_T^{i,n}- \check K_s^{i,n}-\int_s^T \check Z_r^{i,n}dB_r -\int_s^T \int_E \check V_r^{i,n}(e) \tilde{\mu}(dr,de);\\\vspace{0.2cm}
\check Y_s^{i,n} \geqslant  \displaystyle \max_{j \in \mathcal{I}^{-i}}(\check Y_s^{j,n} -\check g_{ij}(s,X_s^{t,x})); \\\vspace{0.2cm}
\textstyle \int_0^T (\check Y_s^{i,n} -\displaystyle \max_{j\in \mathcal{I}^{-i}}(\check Y_s^{j,n} -\check g_{ij}(s,X_s^{t,x}))) d\check K_s^{i,n} = 0.
\end{cases}
 \end{split}
\end{equation}
Next we need to have a representation for
$\check Y^{i,n}$ and $\check V^{i,n}$ similar to \eqref{rep1}, i.e., there exist deterministic continuous bounded functions $(\check u^i)_{\ii}$ such that for any $\tx$ and any $s \in [t,T]$:
 \begin{align}\label{rep1x}
     &\mbox{(a) }\,\check Y_s^{i,n} := \check u^{i,n}(s,\xtx_s)\, \mbox{and }\nonumber\\[4pt]
     & \mbox{(b) }\, \check V_s^{i,n}(e) := \check u^{i, n}(s,\xtx_{s{-}}+ \beta(\xtx_{s{-}},e))- \check u^{i,n}(s,\xtx_{s{-}}),\,ds\otimes d\mathbb{P} \otimes d\lambda \mbox{ on } [t,T] \times \Omega \times E.\nonumber
\end{align}But this can be shown by induction.
 For $n=1$ the property holds true since 
 $\check V_s^{i,0}=0$ and then \\
 $(\check{Y}^{i,1,t,x} ,\check{Z}^{i,1,t,x} ,\check {V}^{i,1,t,x},\check{K}^{i,1,t,x} )_{\ii}$ verify: $\forall s\le T$,
 \begin{equation}\label{recurrencex2n=1}
\begin{split}
\begin{cases}
\vspace{0.2cm} \check Y^{i,1}\in \ss,  \check Z^{i,1} \in \hdd, \check V^{i,1} \in \hld, \mbox{ and } \check K^{i,1}\in \aa;\\ \vspace{0.2cm}
\check Y^{i,1}_s = \check h_i(X_T^{t,x}) +\int_s^T \check  F_i(r,\xtx_r,(\check Y^{k,1}_r)_{k\in \mi},\check Z^{i,1}_r,0)dr \\ \vspace{0.4cm}
 \qquad \qquad \quad \,  +\check K_T^{i,1}- \check K_s^{i,1}-\int_s^T \check Z_r^{i,1}dB_r -\int_s^T \int_E \check V_r^{i,1}(e) \tilde{\mu}(dr,de);\\\vspace{0.2cm}
\check Y_s^{i,1} \geqslant  \displaystyle \max_{j \in \mathcal{I}^{-i}}(\check Y_s^{j,1} -\check g_{ij}(s,X_s^{t,x})); \\\vspace{0.2cm}
\textstyle \int_0^T (\check Y_s^{i,1} -\displaystyle \max_{j\in \mathcal{I}^{-i}}(\check Y_s^{j,1} -\check g_{ij}(s,X_s^{t,x}))) d\check K_s^{i,1} = 0.
\end{cases}
 \end{split}
\end{equation}
Since we are in the Markovian framework and  the functions 
$(\check{F}_i (s,x,\vec{y}, z,0))_{\ii}$ verify (H1)-ii), iii), iv) and \eqref{newcondition} (see Remark \ref{rmqimportante}), and $\check{F}_i(t,x,\vec 0,0,0)$ and $\check{h}_i$, $\ii$, are bounded, then by Theorem \ref{thmexistence2} there exist  deterministic continuous bounded functions $(\check u^{i,1})_{\ii}$ such that for any $\tx$ and $\stt$, $\check Y_s^{i,1} := \check u^{i,1}(s,\xtx_s)$, $\ii$. Next by continuity of $\check u^{i,1}$ and since the L\'evy measure $\lambda(.)$ is finite, we have:
\begin{align*}
    \check V_s^{i,1}(e) := 1_{\{ s\ge t\}}(\check u^{i, 1}&(s,\xtx_{s{-}}+ \beta(\xtx_{s{-}},e))- \check u^{i,1}(s,\xtx_{s{-}})),\,\,ds\otimes d\mathbb{P} \otimes d\lambda \mbox{ on } [0,T] \times \Omega \times E.
\end{align*}
 Therefore the property holds for $n=1$. Next assume that it is satisfied for some $n$. First let us set for any $\ii$,
 \def \bfcek {\check{\underline {F}}}
 $$
 \bfcek_i^n(s,x,\vec y,z):=\check F_i(s,x,\vec y,z,v)_{|v=((\check u^{i,n}(s,x+\beta(x,e)-\check u^{i,n}(s,x))_{e\in E}}. 
 $$
 The generators $\bfcek_i^n(s,x,\vec y,z)$ verify (H1)-ii), iii), iv) and \eqref{newcondition}, moreover 
 $(\bfcek_i^n(s,x,0,0))_{\ii}$ are bounded functions. 
 Therefore by Theorem \eqref{thmexistence2}, there exist bounded deterministic continuous functions $\check u^{n+1}(t,x)$, $\ii$, such that for any $s\in [t,T]$, $\check Y^{i,n+1}_s=\check u^{n+1}(s,\xtx_s)$. Next by continuity of $\check u^{i,n+1}$ and since $\lambda(.)$ is finite, we have for any $\ii$, 
 $
 \check V_s^{i,n+1}(e) := 1_{\{ s\ge t\}}(\check u^{i,n+ 1}(s,\xtx_{s{-}}+ \beta(\xtx_{s{-}},e))- \check u^{i,n+1}(s,\xtx_{s{-}}))
 $. Thus the property holds true for $n+1$ and then (a)-(b) above are satisfied for any $n\ge 0$. 
 
 \noindent To proceed it is enough to follow the same steps as in Steps 3 and 4 in the proof of Theorem \ref{thmexistence2} to show that:

i) Let $\ii$ be fixed. The sequence $(\check u^{i,n})_{n\ge 0}$ converges uniformly on $\TR$ to some bounded continuous function $\check u^i$. The representation given in point (a) above allows to show that the sequence 
$(\check Y^{i,n})_{n\ge 0}$ converges to some process $\bar Y^i$ in $\ss_{[t,T]}$. Next as in Step 4, we have also the convergence of $(\check Y^{i,n})_{n\ge 0}$ in 
$\ss_{[0,t]}$ to $\bar Y^i$ a deterministic continuous bounded function.
Therefore, the sequence $(\check Y^{i,n})_{n\ge 0}$ converges in $\ss$ to some process 
$\bar Y^i$. On the other hand we have also the convergence of $(\check V^{i,n})_{n\ge 0}$ in 
$\hld$ to \\ $\bar V^i(t,e):=
1_{\{ s\ge t\}}(\check u^{i}(s,\xtx_{s{-}}+ \beta(\xtx_{s{-}},e))- \check u^{i}(s,\xtx_{s{-}})),\,\, ds\otimes d\mathbb{P} \otimes d\lambda \mbox{ on } [0,T] \times \Omega \times E$ and the convergence of 
$(\check Z^{i,n})_{n\ge 0}$ (resp. $(\check K^{i,n})_{n\ge 0}$) in $\hdd$ (resp. $\ss$) to a process $\bar Z^{i}$ (resp. $\bar K^{i}$), $\ii$; 

ii) $((\overline{Y}^{i,t,x} ,\overline{Z}^{i,t,x} ,\overline{V}^{i,t,x},\overline{K}^{i,t,x} ))_{\ii}$ is a solution of the system associated with $((\check{F_i})_{\ii}, (\check{h}_i)_{\ii},\check{g}_{ij})_{i,j\in \mathcal{I}})$. 

\noindent To proceed for $s\in [0,T]$, let us set:
 \begin{align*}
     Y_s^{i,t,x}&:= (\phi(\xts_s))^{-1} \overline{Y}_s^{i,t,x},\\[3pt]
      dK_s^{i,t,x}&:= (\phi(\xts_{s-}))^{-1} d\overline{K}_s^{i,t,x} \mbox{ and }K_0^{i,t,x}=0,\\[3pt]
     Z_s^{i,t,x}&:= (\phi(\xts_{s}))^{-1}\big\lbrace \overline{Z}_s^{i,t,x} - ((\phi(\xts_{s}))^{-1} \overline{Y}_s^{i,t,x}D_x\phi(x) \sigma(s,x)\big\rbrace,\\[3pt]
     V_s^{i,t,x}(e) &:= (\phi(X_{s-}^{t,x}+ \beta(\xtx_{s{-}},e))^{-1}\big\lbrace \overline{V}_s^{i,t,x} - \phi(X_{s-}^{t,x})^{-1}\overline{Y}_s^{i,t,x} \big(\phi(X_{s-}^{t,x}+ \beta(\xtx_{s{-}},e)) -\phi(X_{s-}^{t,x})\big)\big\rbrace.
 \end{align*}
Then $(Y^{i,t,x},Z^{i,t,x},V^{i,t,x},K^{i,t,x})_{i\in \mi}$ is a solution of system \eqref{eqBSDE}.
Moreover in setting $u^i(t,x):= (\phi(x))^{-1}\overline{u}^i(t,x),\, (t,x)\in [0,T]\times \R^k$  and $\ii$ we obtain that 
for any $\stt$, $Y^{i,t,x}_s=u^i(s,\xtx_s)$ for any $\ii$ and $u^i$ is of polynomial growth as 
$\bar u^i$ is bounded. $\Box$

As a by-product of the Proposition 3.2 and Theorem 3.3 we have the following:
\begin{corollaire}\label{coro1}
For any $\ii$ and $(t,x)\in [0,T]\times \R^k$,
\begin{align*}
V_s^{i,t,x}(e) =  1_{\{s\ge t\}}\{u^{i}(s,\xtx_{s{-}}&+ \beta(\xtx_{s{-}},e))- u^{i}(s,\xtx_{s{-}})\},\, ds\otimes d\mathbb{P} \otimes d\lambda \mbox{ on } [0,T] \times \Omega \times E.
\end{align*}
\end{corollaire}
Now, we provide the uniqueness of the Markovian solution of the system of reflected BSDEs \eqref{eqBSDE}. 
\begin{proposition}\label{unicité} 
Let $(\tilde{u}^i)_{\ii}$ be   deterministic continuous functions of polynomial growth such that 
\begin{equation}
 \forall s \in [t,T], \, \, Y_s^{i,t,x} = \tilde{u}^i(s,X_s^{t,x}).
\end{equation} 
Then, for any $\ii$, $\tilde{u}^i = u^i$.
\end{proposition}
\nd \sol{\tbf{Proof}}: In order to show that the Markovian  solution of the system of reflected BSDEs \eqref{eqBSDE} is unique, we suppose that there exist other continuous with polynomial growth functions $(\tilde{u}^i)_{\ii}$ such that:
\begin{equation*}
 \forall s \in [t,T], \, \, \tilde{Y}_s^{i,t,x} = \tilde{u}^i(s,X_s^{t,x}),
\end{equation*}
where $(\tilde{Y}^{i,t,x})_{\ii}$ is the first component of the solution of the following system of RBSDEs with jumps with interconnected obstacles: for any $i \in \mathcal{I}$ and $\stt$,
\begin{equation}\label{eqBSDE2}
\begin{split}
\begin{cases}
\vspace{0.3cm} \tilde{Y}^{i,t,x}\in \ss,  \tilde{Z}^{i,t,x} \in \hdd, \tilde{V}^{i,t,x} \in \hld, \mbox{ and } \tilde{K}^{i,t,x}\in \aa;\\ \vspace{0.3cm}
\tilde{Y}_s^{i,t,x} = h_i(X_T^{t,x})+ \int_s^T \bar f_i(r,X_r^{t,x},(\tilde{Y}_r^{k,t,x})_{k\in \mi},\tilde{Z}_r^{i,t,x}, \int_E \gamma_i(X_r^{t,x},e)\tilde{V}_r^{i,t,x}(e) \lambda(de))dr\\\vspace{0.3cm}
 \qquad \qquad    +\tilde{K}_T^{i,t,x} - \tilde{K}_s^{i,t,x}-\int_s^T \tilde{Z}_r^{i,t,x}dB_r -\int_s^T \int_E \tilde{V}_r^{i,t,x}(e) \tilde{\mu}(dr,de),\\\vspace{0.3cm}
 \tilde{Y}_s^{i,t,x} \geqslant  \displaystyle \max_{j \in \mathcal{I}^{-i}}(\tilde{Y}_s^{j,t,x} -g_{ij}(s,X_s^{t,x})),\\\vspace{0.3cm}
 \textstyle {\int_t^T (\tilde{Y}_s^{i,t,x} -\displaystyle \max_{j\in \mathcal{I}^{-i}}(\tilde{Y}_s^{j,t,x} -g_{ij}(s,X_s^{t,x}))) d\tilde{K}_s^{i,t,x} = 0.}
\end{cases}
\end{split}
\end{equation}
On the other hand, as for any $i \in \lbrace 1,...m\rbrace$, $\tilde{u}^i$ is a continuous function of polynomial growth and since the L\'evy measure $\lambda(.)$ is finite, one has
\begin{equation*}
\tilde{V}_s^{i,t,x}(e) =  \tilde{u}^{i}(s,\xtx_{s-}+ \beta(\xtx_{s-},e))- \tilde{u}^{i}(s,\xtx_{s-}),\,\, ds\otimes d\mathbb{P} \otimes d\lambda \mbox{ on } [t,T] \times \Omega \times E.
\end{equation*}
Now, let $\stt$ and 
an admissible strategy $\d \in \mathcal{A}_s^i$. Let $(P^{\d}_r, N^{\d}_r, Q^{\d}_r)_{r\in [s,T]}$ be the triplet of processes  associated with $\d $ and which solves the following BSDE: $\forall  r\in [s,T]$
  \begin{equation*}
  P_r^{\d} = h^{\d}(X_T^{t,x}) +\textstyle \int_r^T f^{\d}(\t,X_\t^{t,x},N^\d_\t) d\t -\textstyle \int_r^T N_\t^{\d} dB_\t - \textstyle\int_r^T \int_E Q_\t^{\d}(e) \tilde{\mu}(d\t,de) - A_T^{\d} + A_r^{\d},  
  \end{equation*}
  where, when $\d_\t =i$, $f^{\d}(\t,X_\t^{t,x},z)$ is equal 
  to 
  \begin{equation*}
  \bar{f}_i(\t,X_\t^{t,x}, (\tilde{u}^k(\t,X_\t^{t,x}))_{k \in \mi},z, \textstyle \int_E \gamma_i(X_\t^{t,x},e) \{\tilde{u}^{i}(\t,\xtx_{\t-}+ \beta(\xtx_{\t-},e))- \tilde{u}^{i}(\t,\xtx_{\t-})\} \lambda(de)).
  \end{equation*}
  Therefore, we have the following representation of $\tilde{Y}^i$:
\begin{equation*}
\tilde{Y}_s^{i} = \displaystyle \mbox{esssup}_{\d \in \mathcal{A}_s^i}(P_s^{\d} - A_s^{\d}). 
\end{equation*}  
Next, the same procedure as the one which leads to inequality \eqref{estimationu} can be used here to deduce that for any $\ii$,
\begin{equation*}\label{estiu}
\begin{aligned}
|u^{i}(t,x) - &\tilde{u}^{i}(t,x)|^2 \leq  2C  \E \big[ \textstyle \int_t^T e^{\alpha_0 (r-t)}\int_E \displaystyle{\sum_{k=1,m} }\big\{ |(u^{i}- \tilde{u}^{i})(r,\xtx_{r-} + \beta(\xtx_{r-},e))|^2 
+ |(u^{i} -\tilde{u}^{i})(r,\xtx_{r-})\big|^2\} \lambda(de)dr\big].
\end{aligned}
\end{equation*}
We now consider two cases. 

\nd  \textbf{\sol{Case 1}:} The functions $u^i$ and $\tilde u^i$, $\ii$, are bounded. 

Let $\eta$ be the constant given in Step 3 and which does not depend on the terminal condition $(h_i)_{\ii}$ and verifies $ \frac{4}{\alpha_0}m\lambda(E)(e^{\alpha_0 \eta }-1)= \frac{3}{4}$. Then, we deduce from \eqref{estimationu}, that for any $\ii$, 
$$
\|u^{i}- \tilde{u}^{i}\|^2_{\infty,\eta}\leq  \frac{3}{4}\|u^{i}- \tilde{u}^{i}\|^2_{\infty,\eta}
$$
which implies that, for any $\ii$, $u^i = \tilde{u}^i$ on $[T-\eta, T]$. Consequently, for any  $s \in [T-\eta,T]$ and $\ii$, \vspace{0.1cm}$Y_s^{i,t,x} = \tilde{Y}_s^{i,t,x}$. Next, on $[T-2\eta, T-\eta]$, we have  
$$
\|u^{i}- \tilde{u}^{i}\|^2_{\infty,2\eta}\leq  \frac{3}{4}\|u^{i}- \tilde{u}^{i}\|^2_{\infty,2\eta}+\frac{4}{\alpha_0}m\lambda(E)(e^{2\alpha_0 \eta}-e^{\alpha_0\eta}) \|u^{i}- \tilde{u}^{i}\|^2_{\infty,\eta}.
$$
Since $u^i = \tilde{u}^i$ on $[T-\eta, T]$, we then obtain:
$$
\|u^{i}- \tilde{u}^{i}\|^2_{\infty,2\eta}\leq  \frac{3}{4}\|u^{i}- \tilde{u}^{i}\|^2_{\infty,2\eta}.$$
Consequently, for any $\ii$, $u^i = \tilde{u}^i$ on $[T-2\eta, T-\eta]$. Thus, for any  $s \in [T-2\eta,T-\eta]$ and $\ii$, \vspace{0.1cm}$Y_s^{i,t,x} = \tilde{Y}_s^{i,t,x}$.
 Repeating now this procedure on $[T-3\eta,T-2\eta]$, $[T-4\eta ,T-3\eta]$ etc., we obtain, for any $\ii$, $u^i = \tilde{u}^i$. Thus, for any $s \in [t,T]$ and $\ii$, $Y_s^{i,t,x} = \tilde{Y}_s^{i,t,x}$. Henceforth, $(Y^{i,t,x})_{\ii}$ is the unique Markovian solution to the system of BSDEs \eqref{eqBSDE}.
 
 \nd  \textbf{\sol{Case 2} :} We deal with the general case, i.e., without assuming the boundedness of the functions 
 $u^i$ and $\tilde u^i$, $\ii$, but only polynomial growth.
 
 Let us define, for $s\in [t,T]$, \begin{equation*}
\overline{Y}_s^{i,t,x} := Y_s^{i,t,x} \phi(X_s^{t,x})\, \,\, \,  \mbox{ and } \, \,\, \,  \underline{Y}_s^{i,t,x} := \tilde{Y}_s^{i,t,x} \phi(X_s^{t,x}),
\end{equation*}
where $\phi$ is the function defined in \eqref{phi}. Therefore 
$(\overline{Y}^{i,t,x},\overline{Z}^{i,t,x},\overline{K}^{i,t,x},
\overline{U}^{i,t,x})_{\ii}$ and \\
$(\underline{Y}^{i,t,x},\underline{Z}^{i,t,x},\underline{K}^{i,t,x},
\underline{U}^{i,t,x})_{\ii}$ are solutions of the system 
\eqref{eqBSDE2} associated with 
$(\check F_i)_{i\in \mathcal{I}}, (\check g_{ij})_{i,j \in \mathcal{I}}$ and $(\check h_i)_{i \in \mathcal{I}}.$
But for any $\ii$, $\overline{Y}^{i,t,x}$ and 
$\underline{Y}^{i,t,x}$ have representations through deterministic continuous bounded functions $\phi u^i$ and $\phi \tilde u^i$, respectively. Therefore by using the result of Step 1 we deduce that 
$\phi u^i=\phi \tilde u^i$ for any $\ii$ and then 
$\overline{Y}_s^{i,t,x}=\underline{Y}^{i,t,x}_s$ for any $s\in [t,T]$ and $\ii$, which implies that 
${Y}_s^{i,t,x} =\tilde{Y}_s^{i,t,x}$, 
for any $s\in [t,T]$ and $\ii$. Thus the Markovian solution of \eqref{eqBSDE} is unique. 
\qed

\section{The main result : Existence and uniqueness of the solution for the system of IPDEs with interconnected obstacles \eqref{local-operator} }

We now turn to the study of the existence and uniqueness in viscosity sense of the solution of the system of integral-partial differential  equations with interconnected obstacles \eqref{system}. Before doing so, we first precise the meaning of the definition of the viscosity solution of this system. It is not exactly the same as in \cite{hamadene2015viscosity} (see also Definition \eqref{systemdef1} in the Appendix). 
 
\begin{definition}\label{nvdef}
We say that a family of deterministic continuous functions $\vec u:=(u^i)_{i\in \mi}$  is a viscosity supersolution (resp. subsolution) of \eqref{system} if: 
$\forall i \in \lbrace 1,...,m\rbrace$, 
\begin{align*}
& \mbox{a) } \,u^i(T,x) \geq (\mbox{resp.} \leq )\,\, h_i(x) , \,\,\forall x \in \R^k\,\,; \\
& \mbox{b)  if}\,\phi \in {\cal C}^{1,2}([0,T] \times \R^k)\mbox{  is 
such that $(t,x) \in [0,T) \times \R^k$ a global minimum}\\
& \qquad \qquad \qquad \mbox{ (resp. maximum) point of $u^i - \phi$}
\end{align*}
then
\begin{align*}
\min \Big\lbrace &u^i(t,x) - \displaystyle \max_{j \in {\mathcal{I}^{-i}}}(u^j(t,x)-g_{ij}(t,x)) ;-\partial_t\phi(t,x) - \mathcal{L}\phi(t,x) - \mathcal{K}\phi(t,x)\\
 & - {\bar f_i}(t,x,(u^k(t,x))_{k=1,...,m},(\sigma^\top D_x\phi)(t,x), \mathcal{B}_iu^i(t,x))\Big\rbrace \geq \,\,(resp. \le )\,\,0. 
\end{align*}
We say that $\vec u:=(u^i)_{i\in \mi}$ is a viscosity solution of \eqref{system} if it is both  a supersolution and subsolution of \eqref{system}.
\end{definition}

\begin{remarque}
In our definition, the last argument of $\bar f_i(.)$ is $\mathcal{B}_iu^i(t,x)$ instead of  $\mathcal{B}_i\phi(t,x)$, where $\phi$ is the test function. Indeed, $\mathcal{B}_iu^i(t,x)$ is well-posed since $u^i$ has a polynomial growth, $\beta$ is bounded and the measure $\lambda(.)$ is finite. 
\end{remarque}
We are now able to state the main result of this paper.

Let $(Y^{i,t,x},Z^{i,t,x},V^{i,t,x},K^{i,t,x})_{\ii}$ be the solution of \eqref{eqBSDE} and let $(u^i)_{\ii}$ be the continuous functions with polynomial growth such that for any $\tx$, $\ii$ and $s\in [t,T]$, $$Y^{i,t,x}_s=u^i(s,\xtx_s).$$
We then have:
\begin{theoreme} Assume that the functions $(\bar f_i)_{i\in \mathcal{I}}$ and $(\gamma_i)_{\ii}$ verify Assumption (H1) and, $(g_{ij})_{i,j \in \mathcal{I}}$  and $(h_i)_{i \in \mathcal{I}}$ verify Assumptions (H2) and (H3).
Then the functions $(u^i)_{\ii}$ is the unique viscosity solution of the system \eqref{system}, according to Definition \eqref{nvdef}, in the class of continuous functions with polynomial growth. 
\end{theoreme}
\nd \sol{\tbf{Proof}}: We first show that  $(u^i)_{\ii}$ is a viscosity solution of system
\eqref{system}. So let us consider the following system of reflected BSDEs: $\fst$,
\begin{equation}\label{nvRBSDE}
\begin{split}
\begin{cases}
\vspace{0.3cm} \underbar{Y}^{i,t,x}\in \ss,  \underbar{Z}^{i,t,x} \in \hdd, \underbar{V}^{i,t,x} \in \hld, \mbox{ and } \underbar{K}^{i,t,x}\in \aa;\\ \vspace{0.3cm}
\underbar{Y}_s^{i,t,x} = h_i(X_T^{t,x})+ \int_s^T {\bar f_i}(r,X_r^{t,x},(\underbar{Y}_r^{k,t,x})_{k\in \mi},\underbar{Z}_r^{i,t,x}, \int_E \gamma_i(\xtx_r,e)\times \\\vspace{0.3cm}
 \quad \quad \quad  \{u^{i}(r,\xtx_{r-}+ \beta(\xtx_{r-},e)) - u^{i}(r,\xtx_{r-})\} \lambda(de))dr+\underbar{K}_T^{i,t,x} - \underbar{K}_s^{i,t,x} -\int_s^T \underbar{Z}_r^{i,t,x}dB_r -\int_s^T \int_E \underbar{V}_r^{i,t,x}(e) \tilde{\mu}(dr,de);\\\vspace{0.3cm}
 \underbar{Y}_s^{i,t,x} \geqslant  \displaystyle \max_{j \in \mathcal{I}^{-i}}(\underbar{Y}_s^{j,t,x} -g_{ij}(s,X_s^{t,x}));\\\vspace{0.3cm}
 \textstyle {\int_0^T (\underbar{Y}_s^{i,t,x} -\displaystyle \max_{j\in \mathcal{I}^{-i}}(\underbar{Y}_s^{j,t,x} -g_{ij}(s,X_s^{t,x}))) d\underbar{K}_s^{i,t,x} = 0.}
\end{cases}
\end{split}
\end{equation}
As the deterministic functions $(u^i)_{i\in \mi}$ are continuous and of polynomial growth, $\beta(x,e)$ and $\gamma_i(x,e)$  verify respectively \eqref{eq4} and \eqref{eqgamma} and finally by Theorem \ref{thmexistence1}, the solution of this system exists and is unique, since the functions  $(h_i)_{\ii}$, $(g_{ij})_{i,j\in \mi}$ and 
$$\vspace{0.2cm} (t,x,y,z) \longmapsto {\bar f_i}(t,x,y,z,\textstyle \int_E \gamma_i(x,e)\{u^{i}(t,x+ \beta(x,e))  - u^{i}(t,x)\} \lambda(de))$$
 satisfy the Assumptions (H1)-(H3) and does not depend on $v$. Moreover, again by  Theorem \ref{thmexistence1}, there exist deterministic continuous functions of polynomial growth $(\underbar{u}^{i})_{\ii}$, such that: For any $\ii$ and $s \in [t,T]$,
\begin{equation*}
 \underbar{Y}_s^{i,t,x} = \underbar{u}^{i}(s,\xtx_s).
\end{equation*} 
Finally, using a result by Hamad\`ene-Zhao \cite{hamadene2015viscosity}, we deduce that $(\underbar{u}^{i})_{\ii}$ is a solution in viscosity sense of the following system of  IPDEs with interconnected obstacle:  
\begin{equation}\label{eqnv}
\begin{cases}
\min \lbrace \underbar{u}^i(t,x) - \displaystyle \max_{j \in {\mathcal{I}^{-i}}}(\underbar{u}^j(t,x)-g_{ij}(t,x)) ;-\partial_t\underbar{u}^i(t,x) - \mathcal{L}\underbar{u}^i(t,x)- \mathcal{K}\underbar{u}^i(t,x)\\
 \quad  - {\bar f_i}(t,x,(\underbar{u}^k(t,x))_{k=1,,...,m},(\sigma^\top D_x\underbar{u}^i)(t,x),\mathcal{B}_iu^i(t,x))\rbrace = 0  ;\\
\underbar{u}^i(T,x) = h_i(x),
\end{cases}
\end{equation}
Let us notice that, in this system \eqref{eqnv}, the last component of $\bar {f}_i$ is   $\mathcal{B}_iu^i(t,x)$ and not $\mathcal{B}_i\underbar{u}^i(t,x)$. Next, recall that $(Y^{i,t,x},Z^{i,t,x},V^{i,t,x},K^{i,t,x})_{\ii}$ solves the system of reflected BSDEs with jumps with interconnected obstacles \eqref{eqBSDE}. Therefore, we know, by Corollary \ref{coro1}, that for any $\tx$, $\ii$ and $s\in [t,T]$, 
$$\vspace{0.3cm}
V_s^{i,t,x}(e) =  u^{i}(s,\xtx_{s-}+ \beta(\xtx_{s{-}},e))- u^{i}(s,\xtx_{s{-}}),\,\, ds\otimes d\mathbb{P} \otimes d\lambda \mbox{ on } [t,T] \times \Omega \times E.
$$
Then $(Y^{i,t,x},Z^{i,t,x},V^{i,t,x},K^{i,t,x})_{\ii}$ verify: for any $s\in [t,T]$ and $\ii$,
 \begin{equation}
\begin{split}
\begin{cases}
\vspace{0.35cm}
Y_s^{i,t,x} = h_i(X_T^{t,x}) +\int_s^T {\bar f_i}(r,X_r^{t,x},(Y_r^{k,t,x})_{k \in \mi},Z_r^{i,t,x},\int_E \gamma_i(\xtx_r,e)\times\\\vspace{0.35cm}
 \quad \quad \quad   \{u^{i}(r,\xtx_{r-}+ \beta(\xtx_{r-},e)) - u^{i}(r,\xtx_{r-})\}\lambda(de))dr+K_T^{i,t,x}- K_s^{i,t,x}- \int_s^T Z_r^{i,t,x}dB_r -\int_s^T \int_E V_r^{i,t,x}(e) \tilde{\mu}(dr,de);\\\vspace{0.35cm}
Y_s^{i,t,x} \geqslant  \displaystyle \max_{j \in \mathcal{I}^{-i}}(Y_s^{,t,xj} -g_{ij}(s,X_s^{t,x}));\\\vspace{0.2cm}
\int_0^T (Y_s^{i,t,x} -\displaystyle \max_{j\in \mathcal{I}^{-i}}(Y_s^{j,t,x} -g_{ij}(s,X_s^{t,x}))) dK_s^{i,t,x} = 0.
\end{cases}
\end{split}
\end{equation}
Therefore, by uniqueness of the Markovian solution of the system of reflected BSDEs  \eqref{nvRBSDE}, we deduce that for any $s \in [t,T]$ and $\ii$, $\underbar{Y}^{i,t,x}_s = Y^{i,t,x}_s$. Then, for any  $\ii$, $ \underbar{u}^i = u^i$. Consequently, $(u^i)_{i \in \mi}$ is a viscosity solution of \eqref{system} according to Definition \ref{nvdef}. 

Now, let us show that $(u^i)_{i \in \mi}$ is the unique solution in the class of continuous functions with polynomial growth. It is based on the uniqueness of the Markovian  solution of the system of reflected BSDEs \eqref{eqBSDE}.

So let $(\bar u^i)_{i \in \ii}$ be another continuous with polynomial growth solution  of \eqref{system} in the sense of Definition \ref{nvdef}, i.e.,
for any $i \in \mathcal{I}$,
\begin{equation}\label{eqnv2}
\begin{cases}
\min \lbrace \bar{u}^i(t,x) - \displaystyle \max_{j \in {\mathcal{I}^{-i}}}(\bar{u}^j(t,x)-g_{ij}(t,x)) ;-\partial_t\bar{u}^i(t,x) - \mathcal{L}\bar{u}^i(t,x)- \mathcal{K}\bar{u}^i(t,x)\\
 \quad  - {\bar f_i}(t,x,(\bar{u}^k(t,x))_{k=1,,...,m},(\sigma^\top D_x\bar{u}^i)(t,x),\mathcal{B}_i\bar u^i(t,x))\rbrace = 0  ;\\
\bar{u}^i(T,x) = h_i(x).
\end{cases}
\end{equation}
Next, let us consider the following system of reflected BSDEs:
\begin{equation}\label{nvRBSDE2}
\begin{split}
\begin{cases}
\vspace{0.3cm} \bar{Y}^{i,t,x}\in \ss,  \bar{Z}^{i,t,x} \in \hdd, \bar{V}^{i,t,x} \in \hld, \mbox{ and } \bar{K}^{i,t,x}\in \aa;\\ \vspace{0.3cm}
\bar{Y}_s^{i,t,x} = h_i(X_T^{t,x})+ \int_s^T {\bar f_i}(r,X_r^{t,x},(\bar{Y}_r^{k,t,x})_{k\in \mi},\bar{Z}_r^{i,t,x}, \int_E \gamma_i(\xtx_r,e)\times \\\vspace{0.3cm}
 \quad\quad \quad   \{\bar u^{i}(r,\xtx_{r-}+ \beta(\xtx_{r-},e)) - \bar u^{i}(r,\xtx_{r-})\} \lambda(de))dr+\bar{K}_T^{i,t,x} - \bar{K}_s^{i,t,x} -\int_s^T \bar{Z}_r^{i,t,x}dB_r -\int_s^T \int_E \bar{V}_r^{i,t,x}(e) \tilde{\mu}(dr,de);\\\vspace{0.3cm}
 \bar{Y}_s^{i,t,x} \geqslant  \displaystyle \max_{j \in \mathcal{I}^{-i}}(\bar{Y}_s^{j,t,x} -g_{ij}(s,X_s^{t,x}));\\\vspace{0.3cm}
 \textstyle {\int_t^T (\bar{Y}_s^{i,t,x} -\displaystyle \max_{j\in \mathcal{I}^{-i}}(\bar{Y}_s^{j,t,x} -g_{ij}(s,X_s^{t,x}))) d\bar{K}_s^{i,t,x} = 0.}
\end{cases}
\end{split}
\end{equation}
As for the reflected BSDEs \eqref{nvRBSDE}, the solution of the system \eqref{nvRBSDE2} exists and is unique since the deterministic functions $(\bar{u}^i)_{i \in \ii}$ are continuous and of polynomial growth.  Moreover, there exists a deterministic continuous functions of polynomial growth $(v^{i})_{\ii}$, such that:
\begin{equation*}
 \forall s \in [t,T], \, \,  \bar{Y}_s^{i,t,x} = v^{i}(s,\xtx_s).
\end{equation*} 
and
\begin{equation}\label{repVbar}
\bar{V}_s^{i,t,x}(e) =  v^{i}(s,\xtx_{s{-}}+ \beta(\xtx_{s{-}},e))- v^{i}(s,\xtx_{s{-}}),\,\, ds\otimes d\mathbb{P} \otimes d\lambda \mbox{ on } [t,T] \times \Omega \times E.
\end{equation}
Then, by using a result by  Hamad\`ene-Zhao \cite{hamadene2015viscosity}, $(v^{i})_{\ii}$ is the unique viscosity solution, in the class of continuous functions with polynomial growth, of the following system:  
\begin{equation}\label{eqnv3}
\begin{cases}
\min \lbrace v^i(t,x) - \displaystyle \max_{j \in {\mathcal{I}^{-i}}}(v^j(t,x)-g_{ij}(t,x)) ; -\partial_tv^i(t,x) - \mathcal{L}v^i(t,x)- \mathcal{K}v^i(t,x)\\
 \quad - {\bar f_i}(t,x,(v^k(t,x))_{k=1,,...,m},(\sigma^\top D_xv^i)(t,x),\mathcal{B}_i\bar u^i(t,x))\rbrace = 0  ;\\
v^i(T,x) = h_i(x),
\end{cases}
\end{equation}
Now, as the functions $(\bar{u}^i)_{i \in \ii}$ solve system \eqref{eqnv3}, hence by uniqueness of the solution of this system \eqref{eqnv3} (see \cite{hamadene2015viscosity}, Proposition 4.2), for any $\ii$ one deduces that
$ \bar{u}^i = v^i $. Next, by the characterization  of the jumps \eqref{repVbar}, for any $\ii$, it holds:
\begin{equation}
\bar{V}_s^{i,t,x}(e) =  \bar{u}^{i}(s,\xtx_{s{-}}+ \beta(\xtx_{s{-}},e))- \bar{u}^{i}(s,\xtx_{s{-}}),\,\, ds\otimes d\mathbb{P} \otimes d\lambda \mbox{ on } [t,T] \times \Omega \times E.
\end{equation}
Going back now to \eqref{nvRBSDE2} and replace the quantity $\bar{u}^{i}(s,\xtx_{s{-}}+ \beta(\xtx_{s{-}},e))- \bar{u}^{i}(s,\xtx_{s{-}})$ with $\bar{V}_s^{i,t,x}(e)$, yields: for any $\ii$ and $s \in [t,T]$,
\begin{equation} \label{nvBSDE3}
\begin{split}
\begin{cases}
\vspace{0.3cm} 
\bar{Y}_s^{i,t,x} = h_i(X_T^{t,x})+ \int_s^T {\bar f_i}(r,X_r^{t,x},(\bar{Y}_r^{k,t,x})_{k\in \mi},\bar{Z}_r^{i,t,x}, \int_E \gamma_i(\xtx_r,e)\bar{V}_r^{i,t,x}(e)\lambda(de))dr \\\vspace{0.3cm}
 \qquad \quad \quad    +\bar{K}_T^{i,t,x} - \bar{K}_s^{i,t,x} -\int_s^T \bar{Z}_r^{i,t,x}dB_r -\int_s^T \int_E \bar{V}_r^{i,t,x}(e) \tilde{\mu}(dr,de);\\\vspace{0.3cm}
 \bar{Y}_s^{i,t,x} \geqslant  \displaystyle \max_{j \in \mathcal{I}^{-i}}(\bar{Y}_s^{j,t,x} -g_{ij}(s,X_s^{t,x}));\\\vspace{0.3cm}
 \textstyle {\int_t^T (\bar{Y}_s^{i,t,x} -\displaystyle \max_{j\in \mathcal{I}^{-i}}(\bar{Y}_s^{j,t,x} -g_{ij}(s,X_s^{t,x}))) d\bar{K}_s^{i,t,x} = 0.}
\end{cases}
\end{split}
\end{equation}
But $(Y^{i,t,x},Z^{i,t,x},K^{i,t,x},V^{i,t,x})_{\ii}$ is a solution of system \eqref{nvBSDE3} and 
$Y^{i,t,x}$ is Markovian. Then, by the uniqueness result of Proposition \eqref{unicité}, one deduces that
\begin{equation*}
\forall \ii,\, \, \bar Y^{i,t,x}_s = Y^{i,t,x}_s, \,\,\forall s\in [t,T]. 
\end{equation*} 
Hence, for any $i \in \mathcal{I}$ and
$\tx$, $ Y^{i,t,x}_t=\bar Y^{i,t,x}_t= u^i(t,x) = \bar{u}^i(t,x) = v^i(t,x)$ which means that the solution of \eqref{system}, according to Definition \eqref{nvdef}, is unique in the class of continuous functions with polynomial growth. \qed

\noindent \textbf{Appendix}

In the paper by Hamad\`ene and Zhao \cite{hamadene2015viscosity}, the definition of the viscosity solution of the system \eqref{system}, is given as follows:  
\begin{definition}\label{systemdef1}
Let $\vec u:=(u^i)_{i\in \mi}$ be a function of $C([0,T] \times \R^k ; \R^m)$.\\
(i)  We say that $\vec u$ is a viscosity supersolution (resp. subsolution) of \eqref{system} if: 
$\forall i \in \lbrace 1,...,m\rbrace$, 
\begin{align*}
& \mbox{a) } \,u^i(T,x) \geq (\mbox{resp.} \leq )\,\, h_i(x) , \,\,\forall x \in \R^k\,\,; \\
& \mbox{b)  if}\,\phi \in {\cal C}^{1,2}([0,T] \times \R^k)\mbox{  is 
such that $(t,x) \in [0,T) \times \R^k$ a global minimum} \\
& \qquad \quad \mbox{ (resp. maximum) point of $u^i - \phi$},
\end{align*}
then
\begin{align*}
\min \Big\lbrace &u^i(t,x) - \displaystyle \max_{j \in {\mathcal{I}^{-i}}}(u^j(t,x)-g_{ij}(t,x)) ;-\partial_t\phi(t,x) - \mathcal{L}\phi(t,x) - \mathcal{K}\phi(t,x)\\
 &- \bar f_i(t,x,(u^k(t,x))_{k=1,...,m},(\sigma^\top D_x\phi)(t,x), \mathcal{B}_i\phi(t,x))\Big\rbrace \geq \,\,(resp. \le )\,\,0. 
\end{align*}
(ii) We say that $\vec u:=(u^i)_{i\in \mi}$ is a viscosity solution of \eqref{system} if it is both  a supersolution and subsolution of \eqref{system}.
\end{definition}

\end{document}